\documentclass{article}

\usepackage{arxiv}
\usepackage[applemac]{inputenc} 		
\usepackage[T1]{fontenc}    		
\usepackage[colorlinks]{hyperref}      
\usepackage{url}            			
\usepackage{amsthm}
\usepackage{booktabs}       		
\usepackage{amsfonts}       		
\usepackage{amsmath}
\usepackage{nicefrac}       		
\usepackage{microtype}      		
\usepackage{lipsum}				
\usepackage[square,numbers]{natbib}
\usepackage{mathtools}
\usepackage{algpseudocode}
\usepackage{graphicx}
\usepackage{indentfirst,latexsym,bm}
\usepackage{amsmath}
\usepackage{amssymb}
\usepackage{xcolor}
\usepackage{comment}
\usepackage{enumitem}
\usepackage{dsfont}
\usepackage{amssymb}

\usepackage{algorithm}
\usepackage{algorithmicx}
\usepackage{float}
\usepackage{subfig}
\usepackage{caption}
\usepackage{enumitem}
\usepackage{todonotes}
\usepackage{makecell}
\usepackage[flushleft]{threeparttable} 
\usepackage{bbding}
\usepackage{tensor}

\usepackage{tikz}
\usetikzlibrary{arrows.meta,positioning}
\usepackage{pgfplots}
\pgfplotsset{compat=newest}
\usepgfplotslibrary{fillbetween}
\usetikzlibrary{shapes,decorations}
\usetikzlibrary{fit}

\colorlet{color1}{blue}
\colorlet{color2}{red!50!black}

\definecolor{ivory}{RGB}{218,215,203}

\definecolor{cuhkp}{RGB}{98,56,105} 	
\definecolor{cuhkpl}{RGB}{152,24,147} 	
\definecolor{cuhkb}{RGB}{219,160,1} 	
\definecolor{cuhkbd}{RGB}{178,129,0} 	
\definecolor{cuhkr}{RGB}{88,35,155}  	

\usepackage{hyperref}[6.83]
\hypersetup{
  colorlinks=true,
  frenchlinks=false,
  pdfborder={0 0 0},
  naturalnames=false,
  hypertexnames=false,
  breaklinks,
  allcolors = cuhkbd,
  urlcolor = color2,
}

\RequirePackage[capitalize,nameinlink]{cleveref}[0.19]

\crefname{section}{section}{sections}
\crefname{subsection}{subsection}{subsections}

\Crefname{figure}{Figure}{Figures}

\crefformat{equation}{\textup{#2(#1)#3}}
\crefrangeformat{equation}{\textup{#3(#1)#4--#5(#2)#6}}
\crefmultiformat{equation}{\textup{#2(#1)#3}}{ and \textup{#2(#1)#3}}
{, \textup{#2(#1)#3}}{, and \textup{#2(#1)#3}}
\crefrangemultiformat{equation}{\textup{#3(#1)#4--#5(#2)#6}}%
{ and \textup{#3(#1)#4--#5(#2)#6}}{, \textup{#3(#1)#4--#5(#2)#6}}{, and \textup{#3(#1)#4--#5(#2)#6}}

\Crefformat{equation}{#2Equation~\textup{(#1)}#3}
\Crefrangeformat{equation}{Equations~\textup{#3(#1)#4--#5(#2)#6}}
\Crefmultiformat{equation}{Equations~\textup{#2(#1)#3}}{ and \textup{#2(#1)#3}}
{, \textup{#2(#1)#3}}{, and \textup{#2(#1)#3}}
\Crefrangemultiformat{equation}{Equations~\textup{#3(#1)#4--#5(#2)#6}}%
{ and \textup{#3(#1)#4--#5(#2)#6}}{, \textup{#3(#1)#4--#5(#2)#6}}{, and \textup{#3(#1)#4--#5(#2)#6}}

\crefdefaultlabelformat{#2\textup{#1}#3}

\newcommand{\mx}{\mathbf{x}}
\newcommand{\mR}{\mathbb{R}}
\newcommand{\my}{\mathbf{y}}
\newcommand{\mY}{\mathcal{Y}}
\newcommand{\mE}{\mathbb{E}}

\newcommand{\mone}{\mathbf{1}}

\newcommand{\bg}{\bar{g}}
\newcommand{\ubg}{\uline{\bar{\mathbf{g}}}}
\newcommand{\bmx}{\bar{\mx}}
\newcommand{\ubmx}{\uline{\bmx}}

\newcommand{\prt}[1]{\left(#1\right)}
\newcommand{\brk}[1]{\left[#1\right]}
\newcommand{\crk}[1]{\left\{#1\right\}}
\newcommand{\norm}[1]{\left\|#1\right\|}
\newcommand{\T}{\intercal}
\newcommand{\mg}{\mathbf{g}}
\newcommand{\inpro}[1]{\left\langle#1\right\rangle}
\newcommand{\proj}{\mathrm{Proj}}
\newcommand{\cL}{\mathcal{L}}

\newcommand{\bsy}[1]{\boldsymbol{#1}}

\newcommand{\bxi}{\bsy{\xi}}
\newcommand{\cF}{\mathcal{F}}
\newcommand{\condE}[2]{\mE\brk{#1 \middle| #2}}

\newcommand{\cC}{\mathcal{C}}

\newcommand{\order}[1]{\mathcal{O}\left(#1\right)}
\newcommand{\orderi}[1]{\mathcal{O}(#1)}
\newcommand{\prti}[1]{(#1)}
\newcommand{\brki}[1]{[#1]}
\newcommand{\crki}[1]{\{#1\}}
\newcommand{\R}{\mathbb{R}}
\newcommand{\condEi}[2]{\mE[#1|#2]}
\newcommand{\normi}[1]{\Vert #1 \Vert}
\newcommand{\sumn}{\sum_{i=1}^n}
\newcommand{\cY}{\mathcal{Y}}
\newcommand{\x}{\mathbf{x}}
\newcommand{\cN}{\mathcal{N}}

\newcommand{\1}{\mathbf{1}}
\newcommand{\md}{\mathbf{d}}

\newcommand{\cT}{\mathcal{T}}
\newcommand{\h}{\mathbf{h}}

\newcommand{\uh}{\underline{h}} 
\usepackage[english]{babel}
\usepackage{amsmath}
\usepackage{graphicx}
\usepackage{float}
\usepackage{caption}

\usepackage{amsfonts}

\usepackage{amsthm}

\usepackage{amsthm}
\numberwithin{equation}{section}

\usepackage{amssymb}

\usepackage{graphicx}
\usepackage{ulem}

\usepackage{tikz}
\usetikzlibrary{matrix,chains,positioning,decorations.pathreplacing,arrows}
\tikzset{%
	every neuron/.style={
		circle,
		draw,
		minimum size=1cm
	},
	neuron missing/.style={
		draw=none, 
		scale=4,
		text height=0.333cm,
		execute at begin node=\color{black}$\vdots$
	},
}

\usepackage{xcolor}
\definecolor{c1}{HTML}{bb7db2}
\definecolor{c2}{HTML}{dee8ec}
\definecolor{c3}{HTML}{fe793d}

\usepackage{enumitem}

\usepackage{pifont}

\newtheorem{theorem}{Theorem}[section]
\newtheorem{corollary}{Corollary}[theorem]
\newtheorem{lemma}[theorem]{Lemma}
\newtheorem{assumption}{Assumption}[section]

\newtheorem{remark}{Remark}[section]

\date{}

\title{Decentralized Min-Max Optimization with Gradient Tracking\thanks{Runze You and Kun Huang contribute equally. This work was partially supported by the National Natural Science Foundation of China under Grant 62003287, and by Shenzhen Science and Technology Program under Grant RCYX20210609103229031.} }

\author{
Runze You\\
The Chinese University of Hong Kong, Shenzhen \\
School of Data Science (SDS) \\
Shenzhen, Guangdong, China \\
\texttt{runzeyou@link.cuhk.edu.cn}
\And
Kun Huang \\
The Chinese University of Hong Kong, Shenzhen \\
School of Data Science (SDS) \\
Shenzhen, Guangdong, China \\
\texttt{kunhuang@link.cuhk.edu.cn}
\And
Shi Pu\\
The Chinese University of Hong Kong, Shenzhen \\
School of Data Science (SDS) \\
Shenzhen, Guangdong, China \\
\texttt{pushi@cuhk.edu.cn}
}

\begin{document}
\maketitle
\begin{abstract}
This paper presents a novel distributed formulation of the min-max optimization problem. Such a formulation enables enhanced flexibility among agents when optimizing their maximization variables. To address the problem, we propose two distributed gradient methods over networks, termed Distributed Gradient Tracking Ascent (DGTA) and Distributed Stochastic Gradient Tracking Ascent (DSGTA). We demonstrate that DGTA achieves an iteration complexity of $\mathcal{O}(\kappa^2\varepsilon^{-2})$, and DSGTA attains a sample complexity of $\mathcal{O}(\kappa^3\varepsilon^{-4})$ for nonconvex strongly concave (NC-SC) objective functions. Both results match those of their centralized counterparts up to constant factors related to the communication network. 
Numerical experiments further demonstrate the superior empirical performance of the proposed algorithms compared to existing methods. 
\end{abstract}

\section{Introduction}
\label{sec:intro}


The min-max optimization problem, $\min_{x\in\mathcal{X}}\max_{y\in\mathcal{Y}} f(x,y)$, is fundamental in various fields such as game theory, operation research, and machine learning. Compared to pure minimization problems, the min-max formulation enables broader applications, including generative adversarial networks (GANs) \cite{goodfellow2014generative}, adversarial learning \cite{madry2017towards}, distributionally robust optimization (DRO) \cite{sinha2017certifying, jin2021non}, and reinforcement learning \cite{dai2018sbeed}. However, solving these problems and achieving satisfactory performance can be challenging due to their large-scale nature. For instance, recent advances in GANs (e.g., GigaGAN \cite{kang2023scaling}) can involve training a billion-scale model on billions of real-world complex Internet images. 

 In particular, this work addresses the distributed min-max optimization problem:
\begin{equation}
	\label{eq:P}
	\min_{x\in\R^d}\max_{\crk{y_i}_{i=1}^n} \frac{1}{n}\sumn f_i(x,y_i), \text{ s.t. }y_i\in\cY_i,
\end{equation}
where $n$ denotes the number of agents, and $f_i:\R^d\times \cY_i\rightarrow \R$ represents the local objective function known by agent $i\in\cN:=\crki{1,2,\ldots,n}$ only.
Compared to most existing works which assume the same variable $y_i$ and the  same constraint set $\cY_i$ for each agent $i$, Problem~\eqref{eq:P} allows agents to have heterogeneous $y_i$ and $\cY_i$, enabling a  different range of applications. 
One typical example is to train a generator $G(x)$ in a distributed GAN model, where each agent $i$ maintains a local real dataset $\mathcal{D}_{\text{real}}^i$ specific to its domain. Then, agent $i$ trains its discriminator $D_{i}(y_i)$ based on its own dataset $\mathcal{D}_{\text{real}}^i$ and the generator $G(x)$. Specifically, the following distributed min-max optimization problem applies:
\begin{equation}
	\label{eq:P_GAN}
	\min_{x\in\R^d}\max_{\crk{y_i\in\R^{d_{i}}}_{i=1}^n} \frac{1}{n}\sumn\crk{\mE_{y_i\sim\mathcal{D}_{\text{real}}^i}\brk{\varphi(D_{i}(y_i))} + \mE_{x\sim\mathcal{D}_x}\crk{\varphi(1 - D_{i}[G(x)])}},
\end{equation}
where $\varphi:\R\rightarrow\R$ is a concave, increasing function, corresponding to different GAN variants (e.g., $\varphi(t) = \log t$ for classical GAN \cite{goodfellow2014generative} or $\varphi(t) = t$ for Wasserstein GAN \cite{arjovsky2017wasserstein}).
Due to the heterogeneity of the local datasets $\mathcal{D}_{\text{real}}^i$ (e.g., varying feature dimensions $d_{i}$),
it is preferable not to enforce consensus on the variables $y_i$. 

Another example pertains to the distributed formulation of the empirical Wasserstein robustness model (WRM) introduced in \cite{sinha2017certifying}:
\begin{equation}
    \label{eq:wrm_intro}
	\begin{aligned}
		\min_{x\in\R^d} \max_{\crk{y_{ij}}_{j=1}^{m_i}\in\cY_i} \frac{1}{n}\sum_{i = 1}^{n} \crk{\frac{1}{m_i}\sum_{j=1}^{m_i} \brk{\ell(x;y_{ij}) - \gamma\norm{y_{ij} - \xi_{ij}}^2_2} },
	\end{aligned}
\end{equation}
where $\crki{\xi_{ij}}_{j=1}^{m_i}$ denote the data samples known by agent $i$ only, {$\ell:\R^d\times \cY_i\rightarrow\R$ is some classifier}, and $\gamma>0$ is a penalty parameter for the $\ell_2$-norm attack. The quantity $\normi{y_{ij}-\xi_{ij}}^2_2$ denotes the ``cost'' for an adversary to perturb $y_{ij}$ to $\xi_{ij}$. In such a case, the data samples and the constraint sets can vary across different agents, thereby leading to different $y_i$'s and $\cY_i$'s. 

To tackle the distributed min-max optimization problem \eqref{eq:P}, we assume in this paper that each function $f_i(x, y)$ is nonconvex in $x$, $\mu$-strongly concave in $y$ (nonconvex strongly concave, NC-SC), and $L$-smooth, with condition number $\kappa:=L/\mu$. 
The goal is to develop algorithms that identify an $\varepsilon$-first-order stationary point ($\varepsilon$-FSP) $\hat{x}$ for the equivalent problem: 
\begin{equation}
	\label{eq:P_min}
	\min_{x\in\R^d} \Phi(x):=\frac{1}{n}\sumn\Phi_i(x),\; \text{with }\Phi_i(x):= \max_{{y_i}\in\cY_i} f_i(x, {y_i}),
\end{equation}
where $\hat{x}$ satisfies $\normi{\nabla \Phi(\hat{x})}^2\leq \varepsilon^2$ (deterministic) or $\mE\brki{\normi{\nabla \Phi(\hat{x})}^2}\leq \varepsilon^2$ (stochastic).

To solve Problem \eqref{eq:P_min} (or equivalently Problem \eqref{eq:P}), we propose two effective decentralized algorithms, termed Distributed Gradient Tracking Ascent (DGTA) and Distributed Stochastic Gradient Tracking Ascent (DSGTA).
Specifically, we demonstrate that DGTA achieves an iteration complexity\footnote{The number of iterations required to attain $\normi{\nabla \Phi(x)}^2\leq \varepsilon^2$.} of $\orderi{\kappa^2\varepsilon^{-2}}$, and DSGTA attains a sample complexity\footnote{The number of stochastic first-order oracle calls necessary to achieve $\mE\brki{\normi{\nabla \Phi(x)}^2}\leq \varepsilon^2$.} of $\orderi{\kappa^3\varepsilon^{-4}}$ for NC-SC objective functions. Both results match those of the centralized methods \cite{lin2020gradient} up to constant factors related to the communication graph; see Table~\ref{tab:comp_full}~and~\ref{tab:comp_s}.

\subsection{Related Works}

The min-max optimization problem has been extensively studied in the literature  under the centralized setting; see e.g., \cite{jin2020local,lin2020gradient,lin2020near,yang2020global,zhang2020single,nouiehed2019solving,yang2022faster,huang2022accelerated} and the references therein. In particular, for nonconvex strongly concave (NC-SC) objective functions, the Gradient Descent Ascent (GDA) method in \cite{lin2020gradient} attains an iteration complexity of $\orderi{\kappa^2\varepsilon^{-2}}$. Such a result was later enhanced to $\orderi{\sqrt{\kappa}\varepsilon^{-2}}$ in \cite{lin2020near}, aligning with the optimal complexity established in \cite{zhang2021complexity}. 
The stochastic variant of GDA (SGDA) attains a sample complexity of $\orderi{\kappa^3 \varepsilon^{-4}}$ \cite{lin2020gradient}. Such a result was improved to $\orderi{\kappa^3\varepsilon^{-3}}$ in \cite{luo2020stochastic} by integrating the variance reduction technique \cite{fang2018spider} with large batch sizes in the order of $\orderi{\varepsilon^{-2}}$. The work in \cite{huang2022accelerated} avoids querying large batches and achieves a sample complexity of $\orderi{\kappa^{4.5}\epsilon^{-3}}$. The dependence on $\varepsilon$ can be further enhanced to reach the lower bound $\orderi{\varepsilon^{-2}}$ given in \cite{zhang2021complexity} for finite sum objective functions \cite{luo2020stochastic}. 

Recently, several distributed algorithms have been proposed to tackle the distributed min-max problem \eqref{eq:P0}  with or without a central coordinator. Notably, such a formulation requires a consensus variable $y$ and a universal constraint set $\cY$ among all the agents, which differ from  Problem \eqref{eq:P}. 
\begin{equation}
	\label{eq:P0}
	\min_{x\in\R^d}\max_{y\in\cY} \frac{1}{n}\sumn f_i(x,y).
\end{equation}
When a central coordinator is available, several algorithms have been considered based on the idea of Federated Averaging (FedAvg) \cite{mcmahan2017communication}, where the coordinator average both $x_i$ and $y_i$ received from the agents. The work in \cite{sharma2022federated} allows for local updates in SGDA and attains a sample complexity of $\orderi{\kappa^4n^{-1}\varepsilon^{-4}}$ for NC-SC objective functions. Such a result was later enhanced to $\orderi{\kappa^3n^{-1}\varepsilon^{-3}}$ under an additional mean squared smoothness (MSS) assumption in \cite{wu2024solving}. The paper \cite{shen2023stochastic} further improves the sample complexity to $\orderi{\kappa^2n^{-1}\varepsilon^{-4}}$ with some smoothing techniques. It is worth noting that the above works all require additional conditions concerning the data heterogeneity, e.g.,  the bounded gradient dissimilarity condition that $\sumn\normi{\nabla f(x,y) - \nabla f_i(x,y)}^2/n\leq \zeta^2$ for some $\zeta^2>0$ \cite{wu2024solving}.

Without a central coordinator,  several decentralized methods have been proposed to address Problem \eqref{eq:P0}. In particular, the methods proposed in \cite{xian2021faster,chen2022simple,zhang2023jointly,mancino2023variance} achieve a sample complexity of $\orderi{\varepsilon^{-3}}$ by combining the gradient tracking method \cite{pu2021distributed,xu2015augmented,di2016next} with various variance reduction techniques \cite{cutkosky2019momentum,fang2018spider}.


The development and study of algorithms that handle Problem \eqref{eq:P} over networked agents are fairly limited and less advanced. When full gradients are available, the GT/DA method in \cite{tsaknakis2020decentralized} attains an iteration complexity of $\tilde{\mathcal{O}}\prti{\varepsilon^{-2}}$ for NC-SC objective functions. However, the method requires multiple inner-loop iterations on $y$, leading to unsatisfactory practical performance. The work in \cite{reisizadeh2020robust} focuses on the robust federated learning problem which can be regarded as a special instance of Problem \eqref{eq:P}. The considered method attains a sample complexity of $\orderi{\varepsilon^{-4}}$ when the objective functions satisfy the nonconvex-Polyak-{\L}ojasiewicz condition and the bounded gradient dissimilarity condition.

\subsection{Main Contribution}
The main contribution of this paper is three-fold. 

Firstly, we consider a novel formulation (Problem~\eqref{eq:P}) for distributed min-max optimization that encompasses a wide range of applications. 
To address such a  problem, we introduce two algorithms, termed Distributed Gradient Tracking Ascent (DGTA) and Distributed Stochastic Gradient Tracking Ascent (DSGTA). Notably, the proposed algorithms do not rely on additional conditions such as the bounded data heterogeneity assumption or the bounded second moments of the stochastic gradients.

Secondly, under nonconvex strongly concave objective functions, we demonstrate that the iteration complexity of DGTA behaves as $\orderi{\kappa^2\varepsilon^{-2}}$, which is comparable to that of  the centralized GDA algorithm in \cite{lin2020gradient} up to constants related to the communication graph (Table~\ref{tab:comp_full}). Under a large batch size, the iteration and sample complexities of DSGTA behave as $\orderi{\kappa^2\varepsilon^{-2}}$ and $\orderi{\kappa^3\varepsilon^{-4}}$, respectively. Notably, both results match those of SGDA in \cite{lin2020gradient} up to constants related to the communication graph (Table~\ref{tab:comp_s}). Moreover, under a constant batch  size $b$, the sample complexity of DSGTA can be improved to $\orderi{\kappa^2\varepsilon^{-4}}$ when $\varepsilon$ is small enough. 
Notably, such a result is independent of the communication graph and outperforms SGDA under a large batch size.

Thirdly, the analysis  of the proposed algorithms allows for different constraint sets $\cY_i$ across the agents, which is  achieved for the first time under the decentralized min-max setting to our knowledge. 
Moreover, we introduce a novel Lyapunov function to handle the coupled recursions, which leads to more tractable and clearer analysis; see Lemma \ref{lem:dsgda_cL} for details.

\begin{table}[htbp]
\centering
\begin{tabular}{@{}cc@{}}
\toprule
Methods                           & Iteration Complexity                                                                                                                                     \\ \midrule
\cite{lin2020gradient} (centralized)          & $\order{\frac{\kappa^2 }{\varepsilon^2}}$                                                                             \\
\cite{tsaknakis2020decentralized} (decentralized) & $\tilde{\mathcal{O}}\prt{\frac{1}{\varepsilon^2}}$                                                                                                 \\
This work (decentralized)                     & $\order{\frac{\kappa^2 }{\sqrt{1-\lambda} \varepsilon^2} + \frac{\kappa }{(1-\lambda)^2\varepsilon^2}}$ \\ \bottomrule
\end{tabular}
\caption{Iteration complexity results of related methods to solve Problem \eqref{eq:P} when full gradients are available. Note that the original result  in \cite{tsaknakis2020decentralized} does not show $\kappa$, $L$ and $(1-\lambda)$ explicitly.
}
	\label{tab:comp_full}
\end{table}

\begin{table}[htbp]
\centering
\setlength{\tabcolsep}{2pt}
\begin{tabular}{@{}cccc@{}}
\toprule
Methods                              & \makecell[c]{Iteration Complexity \\ ($b = \orderi{\kappa\sigma^2\varepsilon^{-2}}$)}                           & \makecell[c]{Sample Complexity \\ ($b = \orderi{\kappa\sigma^2\varepsilon^{-2}}$)}                                         & \makecell[c]{Iteration/Sample Complexity \\ ($b = 1$)}                                                                                                                                                                                             \\ \midrule
\cite{lin2020gradient} (centralized) & $\order{\frac{\kappa^2}{\varepsilon^2}}$                                                       & $\order{\frac{\kappa^3}{\varepsilon^4}}$                                                                  & /                                                                                                                                                                                                                       \\
This work  (decentralized)           & $\order{\frac{\kappa}{(1-\lambda)\varepsilon^2} + \frac{\kappa^2}{(1-\lambda)\varepsilon^2} }$ & $\order{\frac{\kappa^2 }{(1-\lambda)^2 \varepsilon^4} + \frac{\kappa^3}{\sqrt{1-\lambda} \varepsilon^4}}$ & $\mathcal{O}\bigg(\frac{\kappa^2}{\varepsilon^4}\bigg)$ \\ \bottomrule
\end{tabular}
\caption{Sample and iteration complexity results of related methods to solve Problem \eqref{eq:P} using stochastic gradients.  Here $b$ denotes the batch size for computing stochastic gradients. Assume $\varepsilon$ is sufficiently small.}
	\label{tab:comp_s}
\end{table}

\subsection{Notation}
\label{subsec:notation}

Throughout this paper, column vectors are considered by default unless specified otherwise. Let $x_i^t\in\R^d$ denote the iterate of agent $i$ at the $t$-th iteration. We use $\bar{x}\in\R^d$ to denote the averaged variables of $x_i$ among the agents. For instance, $\bar{x}^t:= \sumn x_i^t/n\in\R^d$ represents the average of all the agents' iterates at the $t$-th iteration. For the sake of clarity and presentation, we introduce the stacked variables as follows:
\begin{align*}
	\x^t&:= \brk{x_1^t, x_2^t,\ldots, x_n^t}^{\T}\in\R^{n\times d},\\
	\nabla_x F(\mx^t,\my^t) &:= \left[\nabla_x f_1(x_1^t,y_1^t), \dots, \nabla_x f_n(x_n^t,y_n^t)\right]^{\T} \in \mR^{n\times d},\\
	\h_x^t &:= \left[h_{x,1}^t, h_{x,2}^t,\ldots, h_{x,n}^t\right]^{\T}\in\R^{n\times d},
\end{align*}
which denote the stacked iterates, full gradients (with respect to $x$), and stochastic gradients (with respect to $x$) at the $t$-th iteration, respectively. 

\subsection{Organization}

The rest of this paper is organized as follows. In Section \ref{sec:setup}, we introduce the standing assumptions and the proposed algorithms. We then proceed to conduct preliminary analysis in Section \ref{sec:pre}. The main convergence results for DGTA and DSGTA are presented in Section \ref{sec:main}. Finally, we provide numerical results in Section \ref{sec:sims} and conclude the paper in Section~\ref{sec:conclusion}.

\section{Setup and Algorithms}
\label{sec:setup}
In this section, we introduce the underlying assumptions in subsection \ref{subsec:assumptions} and present the proposed algorithms in subsection \ref{subsec:algs} along with some discussions.

\subsection{Assumptions}
\label{subsec:assumptions}
We consider the standard assumption in the literature regarding the communication network. Let the agents be connected via a graph $\mathcal{G}=(\mathcal{N},\mathcal{E})$ with $\mathcal{N}=\crki{1,2,...,n}$ representing the set of agents and $\mathcal{E}\subseteq \mathcal{N} \times \mathcal{N}$ representing the set of edges connecting the agents. In particular, $(i,i)\in\mathcal{E}$ for all $i\in\mathcal{N}$. The set of neighbors for agent $i$ is denoted by $\mathcal{N}_i=\{j\in \mathcal{N}:(i,j)\in \mathcal{E}\}$.
There is a mixing matrix $W=[w_{ij}]\in\mathbb{R}^{n\times n}$ corresponding to the graph, where $w_{ij}$ represents the weight of the edge $(i,j)$. Specifically, we assume the following conditions for $\mathcal{G}$ and $W$.

\begin{assumption}
	\label{a.network}
	The graph $\mathcal{G}$ is undirected and connected. There exists a link from $i$ and $j$ ($i\neq j$) in $\mathcal{G}$ if and only if $w_{ij} >0$ and $w_{ji}>0$; otherwise, $w_{ij}=w_{ji}=0$. The mixing matrix $W$ is nonnegative, symmetric and stochastic, i.e., $W^{\T} = W$ and $W\mone =\mone$. 
\end{assumption}
Assumption \ref{a.network} guarantees the spectral norm $\lambda$ of the matrix $(W - \mone\mone^{\T}/n)$ is strictly less than one. Then, the spectral gap defined as $(1-\lambda)$ can be used to measure the impact of the network topology on the algorithmic performance. 

Assumptions \ref{a.finite}-\ref{as:set_y} define the properties of the objective functions $f_i(x, y)$ and the constraint set $\cY_i$. In particular, we consider the smooth nonconvex strongly concave (NC-SC) setting. 

\begin{assumption}
	\label{a.finite}
	Each function $f_i(x,y)$ is $L$-smooth, i.e., 
	\begin{align*}
		\norm{\nabla f_i(x, y) - \nabla f_i(x', y')}\leq L \prt{\norm{x - x'} + \norm{y - y'}}
	\end{align*}
for any $(x,y), (x', y') \in \mR^d\times \mY_i$, and  any  $i\in\cN$.

	In addition, each $f_i(x, y)$ is $\mu$-strongly-concave in $y$, i.e., given any $x\in\mR^d$,
	\begin{align*}
		f_i(x, y)\leq f_i(x, y') + \inpro{\nabla_y f_i(x, y'), y-y'} - \frac{\mu}{2}\norm{y - y'}^2
	\end{align*}
 for any $y, y'\in\mY_i$, and any $i\in\cN$.
\end{assumption}

\begin{remark}
	\label{rem:smooth}
	Due to the smoothness of $f_i(\cdot, \cdot)$, we can conclude that $\nabla_x f_i(\cdot, y)$ is also $L$-Lipschitz for any given $y$ since $\|\nabla_x f_i(x,y) - \nabla_x f_i(x',y)\|\le \|\nabla f_i(x,y) - \nabla f_i(x',y)\| \le L\|x - x'\|$. In addition, we have for any $y,y'\in\cY_i$, $x, x'\in\R^d$, and $i\in\cN$ that
	\begin{equation}
		\label{eq:nxL}
		\begin{aligned}
			\norm{\nabla_x f_i(x, y) - \nabla_x f_i(x', y')}^2\leq \norm{\nabla f_i(x, y) - \nabla f_i(x', y')}^2\leq 2L^2\prt{\norm{x - x'}^2 + \norm{y-y'}^2}.
		\end{aligned}
	\end{equation}
	Relation \eqref{eq:nxL} implies that $\nabla_x f_i(\cdot,\cdot)$ is $2L$-Lipschitz continuous. Similarly, we have $\nabla_y f_i(\cdot,\cdot)$ is $2L$-Lipschitz continuous.
\end{remark}
We assume the constraint sets $\cY_i$ is convex and compact for any $i\in\cN$ in Assumption \ref{as:set_y}.
\begin{assumption}
	\label{as:set_y}
	For any $i\in\cN$, the set $\mY_i$ is convex and compact with a diameter $D>0$.
\end{assumption}

It is worth noting that Assumptions \ref{a.finite} and \ref{as:set_y} are standard when considering NC-SC min-max problems; see, e.g., \cite{lin2020gradient,lin2022nonasymptotic}. As a consequence, the functions $\Phi_i(x):= \max_{y\in \mY_i} f_i(x,y)$ are $2\kappa L$-smooth according to \cite[Lemma 4.3]{lin2022nonasymptotic}. This is formally stated in Lemma \ref{lemma12}. 

\begin{lemma}
	\label{lemma12}
	Let Assumptions \ref{a.finite} and \ref{as:set_y} hold. We have that $\Phi_i(\cdot)= \max_{y\in\mY_i}f_i(\cdot, y)$ is $2\kappa L$-smooth with $\nabla \Phi_i(x)= \nabla_x f(x, \hat{y}_i(x))$, where $\hat{y}_i(x):= \arg\max_{y\in\mY_i}f_i(x, y)$. In addition, $\hat{y}_i(\cdot)$ is $\kappa$-Lipschitz, i.e., 
	\begin{align*}
		\norm{\hat{y}_i(x) -  \hat{y}_i(x')} \leq \kappa\norm{x - x'},\ \forall x, x'\in\mR^d, \ i\in\cN.
	\end{align*}
\end{lemma}

The functions $\Phi$ is assumed to be lower bounded as stated in Assumption \ref{a.low}.
\begin{assumption}
	The function $\Phi(x)=\frac{1}{n}\sumn \max_{{ y_i}\in \mY_i} f_i(x,{ y_i})$ is lower bounded, i.e. $\Phi(x) \ge \Phi^* > -\infty$ for any $x \in \mathrm{dom}(\Phi)$. 
	\label{a.low}
\end{assumption}


When the full gradients of the objective functions are not available, we assume that each agent $i$ is able to obtain an estimated gradient $h_i(x,y,\xi_i)$ of $\nabla f_i(x, y)$ given $x\in\mR^p, y\in\mY_i$ by querying a stochastic gradient oracle, where $\xi_i$ represents a random sample drawn by agent $i$.
In particular, we consider the following standard assumption appeared in the literature including \cite{pu2021distributed,lian2017can,huang2023distributed}. The unbiasedness can be ensured by uniform sampling with replacement; see e.g.,~\cite{bottou2018optimization,huang2023distributed}.

\begin{assumption}
	\label{as:sgrad}
	For all $i\in\cN$ and all $x\in\mR^p, y\in\mY_i$, each random vector $\xi_{i}$ is independent (among agents), and 
	\begin{align*}
		&\condE{h_i(x,y,\xi_i)}{(x, y)} = \nabla f_i(x, y),\\
		&\condE{\norm{h_i(x,y,\xi_i) - \nabla f_i(x, y)}^2}{(x, y)}\leq \sigma^2.
	\end{align*}
\end{assumption}
Similar to $h_i(x,y,\xi_i)$, we let $h_{x,i}(x,y,\xi_i)$ and $h_{y,i}(x,y,\xi_i)$ to denote the estimates for $\nabla_x f_i(x, y)$ and $\nabla_y f_i(x, y)$, respectively.
In addition, given a batch of independently drawn samples $\{\xi_{i,j}\}_{j=1}^b$ with size $b$, we define 
\begin{align}
\label{eq:bsgrad}
    \uh_{x, i} & := \frac{1}{b}\sum_{j=1}^b h_{x,i}(x,y,\xi_{i,j}),\text{ and }
	\uh_{y, i} : = \frac{1}{b}\sum_{j=1}^b h_{y,i}(x,y,\xi_{i,j}),
\end{align}
as the mini-batch stochastic gradients of agent $i$.

\subsection{Algorithms}
\label{subsec:algs}

In this part, we first elucidate the motivation behind the proposed algorithms.
Recall that $\Phi_i(\cdot) = \max_{y\in\cY_i} f_i(\cdot, y)$ is a smooth function in light of Lemma \ref{lemma12}, and thus Problem \eqref{eq:P} is equivalent to the smooth minimization problem~\eqref{eq:P_min}.
Consequently, any distributed algorithm can be applied to solving Problem \eqref{eq:P_min} provided that good estimates of the gradients $\nabla \Phi_i(x)= \nabla_x f(x, \hat{y}_i(x))$ (Lemma~\ref{lemma12}) 
are available. 
It has been observed that a single step of gradient ascent on $y_i$ can serve as a reliable substitution to computing $\hat{y}_i(x) = \arg\max_{y\in\cY_i} f_i(x,y)$ \cite{lin2020gradient}. In addition, the gradient tracking method has been proved to be efficient in mitigating the influence of data heterogeneity \cite{pu2021distributed,di2016next,xu2015augmented}. These observations motivate the first algorithm, termed Distributed Gradient Tracking Ascent (DGTA), outlined in Algorithm \ref{DGTA}.

More specifically, each agent in the DGTA method updates its decision variable $x_i$ using the gradient tracking technique in Line \ref{line:dgta_x} and \ref{line:dgta_xgt}. To obtain the estimate for $\hat{y}_i(x_i)$, a projected gradient ascent step is performed in Line \ref{line:dgta_y} and \ref{line:dgta_yproj}. It is worth noting that agent $i$ refrains from performing a consensus update on the variable $y_i$, which differs from previous works such as \cite{xian2021faster,liu2023precision,liu2020decentralized} due to the new formulation \eqref{eq:P}.

\begin{algorithm}
	\caption{Distributed Gradient Tracking Ascent (DGTA)}
	\begin{algorithmic}[1] 
		\Require  Initialize $(x_i^0, y_i^0)$ and set $g_i^0= \nabla_x f_i(x_i^0, y_i^0)$ for each agent $i\in \cN$, determine $W = [w_{ij}]\in\R^{n\times n}$, stepsizes $\eta_x$ and $\eta_y$.
		\For{$t = 0,1,2,\cdots, T-1$}
		\For{Agent $i$ in parallel}
        \State Communicate with its neighbors $j\in\cN_i$ and update $x_i^{t + 1} = \sum_{j\in\cN_i}w_{ij}(x_j^t - \eta_x g_j^t)$. \label{line:dgta_x}
		\State Update $y_i^{t + \frac{1}{2}} = y_i^t + \eta_y\nabla_y f_i(x_i^t, y_i^t)$.\label{line:dgta_y}
		\State Update $y_{i}^{t + 1} =\proj_{\cY_i}(y_i^{t + \frac{1}{2}})$. \label{line:dgta_yproj}
		\State Communicate with its neighbors $j\in\cN_i$ and update 
		$$g_i^{t + 1} = \sum_{j\in\cN_i} w_{ij}(g_j^t + \nabla_x f_i(x_i^{t + 1}, y_i^{t + 1}) - \nabla_x f_i(x_i^{t}, y_i^{t})).$$ \label{line:dgta_xgt}
		\EndFor
		\EndFor
		\Statex \textbf{Output:} $x_i^{T}$ and $y_i^{T}$.
	\end{algorithmic}
	\label{DGTA}
\end{algorithm}

When the full gradients of $f_i$ are not available, we consider a stochastic variant of DGTA and propose the Distributed Stochastic Gradient Tracking Ascent (DSGTA) method outlined in Algorithm \ref{DSGTA}. In this case, at every iteration~$t$, each agent $i$ approximates $\nabla_x f_i(x^t_i,y^t_i)$ and $\nabla_y f_i(x^t_i,y^t_i)$ by drawing a batch of $b$ samples to compute the mini-batch stochastic gradients $\uh_{x,i}^t$ and $\uh_{y,i}^t$ (defined in \eqref{eq:bsgrad}), respectively.

Notably, both algorithms inherit the favorable characteristics of gradient tracking based methods. The gradient tracker $g_i^t$ tracks the averaged gradients $\bar{g}^t = \sumn \nabla_x f_i(x_i^t, y_i^t)/n$ in DGTA and tracks the averaged stochastic gradients $\bar{g}^t = \sumn \uh_{x,i}^t/n$ in DSGTA.

\begin{algorithm}
	\caption{Distributed Stochastic Gradient Tracking Ascent (DSGTA)}
	\begin{algorithmic}[1] 
		\Require   Initialize $(x_i^0, y_i^0)$ for each agent $i\in \cN$, determine $W = [w_{ij}]\in\R^{n\times n}$, mini-batch size $b$, stepsizes $\eta_x$ and $\eta_y$.
  \For{Agent $i$ in parallel}
		\State  Sample $\bxi_i^{0}=\{\xi_{i,j}^0\}_{j=1}^b$ with size $b$.
  \State Calculate $g^{0}_i = \uh_{x,i}^0 = \frac{1}{b}\sum_{j=1}^b h_{x,i}(x^{0}_i,y^{0}_i,\xi_{i,j}^{0})$ and $\uh_{y,i}^0 = \frac{1}{b}\sum_{j=1}^b h_{y,i}(x_i^0,y_i^0,\xi_{i,j}^0)$.
  \EndFor
		\For{$t = 0,1,2,\cdots, T-1$}
		\For{Agent $i$ in parallel}
        \State Communicate with its neighbors $j\in\cN_i$ and update $x_i^{t + 1} = \sum_{j\in\cN_i}w_{ij}(x_j^t - \eta_x g_j^t)$.\label{line:dsgta_x}
		\State Update $y_i^{t + \frac{1}{2}} = y_i^t + \eta_y \uh_{y,i}^t$
		\State Update $y_{i}^{t + 1} =\proj_{\cY_i}(y_i^{t + \frac{1}{2}})$
		\State Sample $\bxi^{t + 1}_i = \crki{\xi_{i,j}^{t + 1}}_{j=1}^b$ with size $b$. Calculate stochastic gradients $\uh_{x,i}^{t + 1}=\frac{1}{b}\sum_{j=1}^b h_{x,i}(x_i^{t + 1}, y_i^{t + 1};\xi_{i,j}^{t + 1})$ and $\uh_{y,i}^{t + 1} = \frac{1}{b}\sum_{j=1}^b h_{y,i}(x_i^{t + 1}, y_i^{t + 1},\xi_{i,j}^{t + 1})$. \label{line:dsgta_xi}
        \State Communicate with its neighbors $j\in\cN_i$ and update 
		$g_i^{t + 1} = \sum_{j\in\cN_i} w_{ij}(g_j^t + \uh_{x,j}^{t + 1} - \uh_{x,j}^t).$ \label{line:gt}
		\EndFor
		\EndFor
		\Statex \textbf{Output:} $x_i^{T}$ and $y_i^{T}$.
	\end{algorithmic}
	\label{DSGTA}
\end{algorithm}

\section{Preliminary Analysis}
\label{sec:pre}

In this section, we study the convergence properties of DSGTA as outlined in Algorithm \ref{DSGTA}. The convergence of DGTA can be viewed as a special case of DSGTA when $\sigma^2=0$. 

For DSGTA, a direct result is that $\bar{x}^{t + 1} = \bar{x}^t - \eta_x\sumn g_i^t/n$ based on Assumption \ref{a.network}. Such a relation leads to the expected ``approximate'' descent property of the $2\kappa L$-smooth function $\Phi(\cdot) = \frac{1}{n}\sum_{i=1}^n\max_{y\in\mY_i} f_i(\cdot, y)$, stated in Lemma \ref{lem:dsgda_descent}.

\begin{lemma}
	\label{lem:dsgda_descent}
	Let Assumptions  \ref{a.network}-\ref{as:sgrad} hold. Set the stepsize $\eta_x$ to satisfy $\eta_x\leq 1/(2\kappa L)$. We have 
	\begin{equation}
		\label{eq:dsgda_descent}
		\begin{aligned}
			\mE\brk{\Phi(\bar{x}^{t + 1})} &\leq \mE\brk{\Phi(\bar{x}^t)} - \frac{\eta_x}{2}\mE\brk{\norm{\nabla\Phi(\bar{x}^t)}^2} -\frac{\eta_x}{2}\prt{1 - 2\kappa L\eta_x}\mE\brk{\norm{\frac{1}{n}\sumn\nabla_x f_i(x_i^t,y_i^t)}^2} + \frac{\eta_x L^2}{n}\mE\brk{\delta^t}\\
			&\quad + \frac{\eta_x L^2}{n}\mE\brk{\norm{\mx^t - \ubmx^t}^2} + \frac{\eta_x^2 L\kappa \sigma^2}{nb},
		\end{aligned}
	\end{equation}
	where $\delta^t = \sum_{i=1}^n\norm{\hat{y}_i(\bar{x}^t) - y_i^t}^2$ with $\hat{y}_i(\cdot) = \arg\max_{y\in\mY_i} f_i(\cdot, y)$ and $\ubmx^t:= \1(\bar{x}^t)^{\T}$.
\end{lemma}

\begin{proof}
	See Appendix \ref{app:dsgda_descent}.
\end{proof}

Lemma \ref{lem:dsgda_descent} serves as a starting point for the subsequent analysis. Specifically, we are inspired to derive the recursions for the error terms $\mE\brki{\delta^t}$ and $\mE\brki{\normi{\mx^t - \ubmx^t}^2}$. 
The recursion for $\mE\brki{\delta^t}$ is stated below.
\begin{lemma}
	\label{lem:dsgda_delta}
	Suppose Assumptions \ref{a.network}-\ref{as:sgrad} hold. Let $\eta_y\leq 1/[8(L+\mu)]$. We have 
	\begin{equation}
        \label{eq:dsgda_delta}
        \begin{aligned}
        		\mE\brk{\delta^{t + 1}} &\leq \prt{1 - \frac{\eta_y\mu}{4}}\mE\brk{\delta^t} + \frac{9L^2\eta_y}{\mu}\mE\brk{\norm{\mx^t - \ubmx^t}^2} + \frac{3\eta_y^2n\sigma^2}{b} + \frac{4\kappa^2\eta_x^2\sigma^2}{\eta_y\mu b}\\
				&\quad + \frac{4n\kappa^2\eta_x^2}{\eta_y\mu}\mE\brk{\norm{\frac{1}{n}\sumn \nabla_x f_i(x_i^t,y_i^t)}^2}.
        	\end{aligned}
    \end{equation}
\end{lemma}

\begin{proof}
	See Appendix \ref{app:dsgda_delta}.
\end{proof}

Lemma \ref{lem:dsgda_cons0} below states the recursion for $\mE\brki{\normi{\mx^t - \ubmx^t}^2}$.

\begin{lemma}
	\label{lem:dsgda_cons0}
	Let Assumption \ref{a.network} hold. We have 
	\begin{equation*}
	    \mE\brk{\|\mx^{t+1} - \ubmx^{t + 1}\|^2} \le \frac{1+\lambda^2}{2} \mE\brk{\|\mx^t - \ubmx^t\|^2} + \frac{2\lambda^2\eta_x^2}{1-\lambda^2} \mE\brk{\| \mg^t - \ubg^t\|^2},
	\end{equation*}
        where $\ubg^t:= \1(\bar{g}^t)^{\T}$.
\end{lemma}

\begin{proof}
	See Appendix \ref{app:dsgda_cons0}.
\end{proof}

Lemma \ref{lem:dsgda_cons0} further guides us to consider the error term $\mE\brki{\normi{\mg^t - \ubg^t}^2}$, whose recursion is stated in Lemma \ref{lem:dsgda_cons0}.
\begin{lemma}
	\label{lem:dsgda_gtc}
	Suppose Assumptions \ref{a.network}-\ref{as:sgrad} hold. Let 
	\begin{align*}
		\eta_x \leq \frac{(1-\lambda^2)^{\frac{3}{2}}}{6\sqrt{2}L}.
	\end{align*}
	We have 
	\begin{align*}
		\mE\brk{\norm{\mg^{t + 1} - \ubg^{t + 1}}^2} &\leq \frac{3+\lambda^2}{4}\mE\brk{\norm{\mg^t - \ubg^t}^2} + \frac{18 L^2}{1-\lambda^2}\mE\brk{\delta^t}+ \prt{2 + \frac{9\eta_y L^2}{\mu}}\frac{9L^2}{1-\lambda^2}\mE\brk{\norm{\mx^t - \ubmx^t}^2}\\
		&\quad + \frac{9\eta_x^2 nL^2}{1-\lambda^2}\prt{1 + \frac{(4 + \eta_y\mu)\kappa^2}{\eta_y\mu}}\mE\brk{\norm{\frac{1}{n}\sumn\nabla_x f_i(x_i^t,y_i^t)}^2}\\
		&\quad + \brk{2 + \frac{9\eta_y^2 L^2}{1-\lambda^2} + \frac{9\eta_x^2 L^2}{n(1-\lambda^2)}\prt{1 + \frac{(4 + \eta_y\mu)\kappa^2}{\eta_y\mu}}}\frac{n\sigma^2}{b}.
	\end{align*}

\end{lemma}

\begin{proof}
	See Appendix \ref{app:dsgda_gtc}.
\end{proof}

In light of Lemmas \ref{lem:dsgda_descent}-\ref{lem:dsgda_gtc}, we are ready to introduce the Lyapunov function $\cL_t$ defined in \eqref{eq:dsgda_lya_full} and present its recursion in Lemma \ref{lem:dsgda_cL}. Specifically,  define
\begin{equation}
	\label{eq:dsgda_lya_full}
	\begin{aligned}
		\cL^t := \Phi(\bar{x}^t) - \Phi^* + \frac{300\kappa^2 L^2\eta_x}{n(1-\lambda^2)}\norm{\mx^t - \ubmx^t}^2 + \frac{8\eta_xL^2}{n\eta_y\mu}\delta^t + \frac{2400\eta_x^3\kappa^2 L^2}{n(1-\lambda^2)^3}\norm{\mg^t - \ubg^t}^2.
	\end{aligned}
\end{equation}

The detailed derivations for the coefficients in \eqref{eq:dsgda_lya_full} are presented in Appendix~\ref{app:dsgda_At}.  The following lemma describes the recursion for $\cL^t$.

\begin{lemma}
	\label{lem:dsgda_cL}
	Suppose Assumptions \ref{a.network}-\ref{as:sgrad} hold. Let 
    \begin{align*}
        \eta_y\leq \frac{\sqrt{1-\lambda^2}}{9(L+\mu)},\; \eta_x\leq \min\crk{\frac{(1-\lambda^2)^2}{120\sqrt{3}\kappa L}, \frac{\eta_y}{8\sqrt{6}\kappa^2}, \frac{(\eta_y\mu)^{1/4}(1-\lambda^2)}{50\kappa L}}.
    \end{align*}
    Then,
    \begin{equation}
		\label{eq:dsgda_cL}
		\begin{aligned}
				\mE\cL^{t + 1}
				&\leq \mE\cL^t - \frac{\eta_x}{2}\mE\brk{\norm{\nabla \Phi(\bar{x}^2)}^2} + \frac{\eta_x^2 L\kappa \sigma^2}{nb} + \frac{24\eta_x\eta_y\kappa L\sigma^2}{b} \\
				&\quad + \brk{\frac{1}{n\eta_y^2\mu^2} + \frac{225}{(1-\lambda^2)^3} + \frac{54\eta_x^2\kappa^2 L^2}{n(1-\lambda^2)\eta_y\mu} }\frac{32\eta_x^3\kappa^2 L^2\sigma^2}{b}.
			\end{aligned}
	\end{equation}
\end{lemma}

\begin{proof}
	See Appendix \ref{app:dsgda_At}.
\end{proof}

\section{Main Results}
\label{sec:main}

In this section, we present the main results which outline the convergence rate of the DSGTA method (Algorithm~\ref{DSGTA}) for solving Problem \eqref{eq:P} with NC-SC objective functions, detailed in Theorem~\ref{thm:dsgda}. 
The convergence result for the DGTA method (Algorithm~\ref{DGTA}) is a direct implication of Theorem \ref{thm:dsgda} and is stated in Corollary \ref{cor:sigma0}. 
In addition, Corollaries~\ref{cor:comp_small_eta}~and~\ref{cor:comp_large_b} illustrate the convergence rate of DSGTA under a constant batch size and a large batch size that depends on the accuracy level $\varepsilon$, respectively. The complexity results are stated accordingly.

\begin{theorem}
	\label{thm:dsgda}
	Suppose Assumptions \ref{a.network}-\ref{as:sgrad} hold, and let 
	\begin{align*}
        \eta_y\leq \frac{\sqrt{1-\lambda^2}}{9(L+\mu)},\; \eta_x\leq \min\crk{\frac{(1-\lambda^2)^2}{120\sqrt{3}\kappa L}, \frac{\eta_y}{8\sqrt{6}\kappa^2}, \frac{(\eta_y\mu)^{1/4}(1-\lambda^2)}{50\kappa L}}.
    \end{align*}
	We have for the DSGTA method that
	\begin{equation}
		\label{eq:dsgda_avgE}
		\begin{aligned}
			&\frac{1}{T}\sum_{t=0}^{T-1}\mE\brk{\norm{\nabla \Phi(\bar{x}^t)}^2}\leq \frac{2\Delta_\Phi }{\eta_x T} + \frac{600\kappa^2 L^2 \norm{\mx^0 - \ubmx^0}^2}{n(1-\lambda^2) T} + \frac{16D^2L^2}{\eta_y\mu T}+ \frac{4800\eta_x^2\kappa^2 L^2 \sigma^2}{b(1-\lambda^2)^3T} + \frac{2\eta_x L\kappa \sigma^2}{nb}\\
			&\quad \frac{4800\eta_x^2\kappa^2 L^2 \sumn\norm{\nabla_x f_i(x_i^0,y_i^0)}^2}{n(1-\lambda^2)^3T} + \frac{48\eta_y\kappa L\sigma^2}{b}+ \brk{\frac{1}{n\eta_y^2\mu^2} + \frac{225}{(1-\lambda^2)^3} + \frac{54\eta_x^2\kappa^2 L^2}{n(1-\lambda^2)\eta_y\mu} }\frac{64\eta_x^2\kappa^2 L^2\sigma^2}{b},
		\end{aligned}
	\end{equation}
        where $\Delta_\Phi:= \Phi(\bar{x}^0) - \Phi^*$.
\end{theorem}

\begin{proof}
	
	Taking the average over $t= 0,1,\ldots, T-1$ on \eqref{eq:dsgda_cL} leads to 
	\begin{equation}
		\label{eq:dsgda_cL_avg1}
		\begin{aligned}
			\frac{1}{T}\sum_{t=0}^{T-1}\mE\brk{\norm{\nabla \Phi(\bar{x}^t)}^2}&\leq \frac{2\mE\cL^0 }{\eta_x T} + \frac{2\eta_x L\kappa \sigma^2}{nb} +  \frac{48\eta_y\kappa L\sigma^2}{b} \\
			&\quad + \brk{\frac{1}{n\eta_y^2\mu^2} + \frac{225}{(1-\lambda^2)^3} + \frac{54\eta_x^2\kappa^2 L^2}{n(1-\lambda^2)\eta_y\mu} }\frac{64\eta_x^2\kappa^2 L^2\sigma^2}{b}.
		\end{aligned}
	\end{equation}

	Recall that $\Delta_\Phi:= \Phi(\bar{x}^0) - \Phi^*$. We have 
	\begin{equation}
		\label{eq:cL0}
		\begin{aligned}
			\mE\cL^0 &= \Delta_\Phi + \frac{300\kappa^2 L^2\eta_x}{n(1-\lambda^2)}\norm{\mx^0 - \ubmx^0}^2 + \frac{8\eta_xL^2}{n\eta_y\mu}\delta^0 + \frac{2400\eta_x^3\kappa^2 L^2}{n(1-\lambda^2)^3}\mE\brk{\norm{\mg^0 - \ubg^0}^2}\\
			&\leq \Delta_\Phi + \frac{300\kappa^2 L^2\eta_x}{n(1-\lambda^2)}\norm{\mx^0 - \ubmx^0}^2 + \frac{8\eta_xL^2 D^2}{\eta_y\mu}+ \frac{2400\eta_x^3\kappa^2 L^2}{n(1-\lambda^2)^3}\brk{\frac{n\sigma^2}{b} + \sumn\norm{\nabla_x f_i(x_i^0,y_i^0)}^2}.
		\end{aligned}
	\end{equation}

	Substituting \eqref{eq:cL0} into \eqref{eq:dsgda_cL_avg1} leads to the desired result.
	
\end{proof}


Corollary \ref{cor:sigma0} states the convergence result of DSGT, i.e., when the full gradients are available ($\sigma^2=0$). 

\begin{corollary}[DGTA]
	\label{cor:sigma0}
	Suppose Assumptions \ref{a.network}-\ref{as:sgrad} hold. Let 
    \begin{align*}
        \eta_y = \frac{\sqrt{1-\lambda^2}}{9(L+\mu)},\; \eta_x = \frac{1}{\gamma},\;
		\gamma = \frac{120\sqrt{3}\kappa L}{(1-\lambda^2)^2} + \frac{72\sqrt{6}\kappa^2(L+\mu)}{\sqrt{1-\lambda^2}} + \frac{110\kappa^{5/4} L}{(1-\lambda^2)^{9/8}}.
    \end{align*}
    We have for the DGTA method (Algorithm \ref{DGTA}) that
	\begin{equation}
		\label{eq:sigma0}
		\begin{aligned}
			\frac{1}{T}\sum_{t=0}^{T-1}\mE\brk{\norm{\nabla \Phi(\bar{x}^t)}^2}&\leq \frac{2\Delta_\Phi }{\eta_x T} + \frac{600\kappa^2 L^2 \norm{\mx^0 - \ubmx^0}^2}{n(1-\lambda^2) T} + \frac{16D^2L^2}{\eta_y\mu T}\\
			&\quad + \frac{4800\eta_x^2\kappa^2 L^2 \sumn\norm{\nabla_x f_i(x_i^0,y_i^0)}^2}{n(1-\lambda^2)^3T}.
		\end{aligned}
	\end{equation} 

    In particular, 
    \begin{equation}
        \label{eq:sigma0_order}
       \begin{aligned}
         \frac{1}{T}\sum_{t=0}^{T-1}\mE\brk{\norm{\nabla \Phi(\bar{x}^t)}^2}&= \order{\frac{\Delta_\Phi\kappa L}{(1-\lambda)^2 T} + \frac{\Delta_\Phi\kappa^2L}{\sqrt{1-\lambda} T} + \frac{\Delta_\Phi\kappa^{5/4} L}{(1-\lambda)^{9/8}T} + \frac{\kappa D^2L^2}{\sqrt{1-\lambda}T}\right.\\
         &\left.\quad + \frac{\kappa^2 L^2 \norm{\mx^0 - \ubmx^0}^2}{n(1-\lambda) T}  + \frac{ \sumn\norm{\nabla_x f_i(x_i^0,y_i^0)}^2}{nT}},
       \end{aligned}
    \end{equation}
    where $\orderi{\cdot}$ hides the numerical constants.

\end{corollary}
\begin{proof}
    Letting $\sigma^2 = 0$ in \eqref{eq:dsgda_avgE} yields \eqref{eq:sigma0}. From the choices of $\eta_x$ and $\eta_y$, we have $\eta_y\mu = \sqrt{1-\lambda^2}/[9(\kappa+1)]$ and 
	\begin{equation}
		\label{eq:eta_inv_sigma0}
		\begin{aligned}
			\frac{\Delta_\Phi}{\eta_x T} &= \frac{120\Delta_\Phi\sqrt{3}\kappa L}{(1-\lambda^2)^2 T} + \frac{72\Delta_\Phi\sqrt{6}\kappa^2(L+\mu)}{\sqrt{1-\lambda^2} T} + \frac{110\Delta_\Phi\kappa^{5/4} L}{(1-\lambda^2)^{9/8}T},\\
			\eta_x^2&\leq \frac{(1-\lambda^2)^4}{14400 \kappa^2 L^2}.
		\end{aligned}
	\end{equation}

    Substituting \eqref{eq:eta_inv_sigma0} into \eqref{eq:sigma0} yields the desired result \eqref{eq:sigma0_order}.
\end{proof}
\begin{remark}
It can be seen from Corollary \ref{cor:sigma0} that the convergence rate of DGTA depends on the spectral gap $1-\lambda$. A larger gap that implies a better connected communication graph leads to faster convergence. 
\end{remark}
\begin{remark}
    From Corollary \ref{cor:sigma0}, the iteration complexity of DGTA behaves as 
    \begin{align*}
        &\order{\frac{\Delta_\Phi\kappa L}{(1-\lambda)^2 \varepsilon^2} + \frac{\Delta_\Phi\kappa^2L}{\sqrt{1-\lambda}  \varepsilon^2} + \frac{\Delta_\Phi\kappa^{5/4} L}{(1-\lambda)^{9/8} \varepsilon^2} + \frac{\kappa D^2L^2}{\sqrt{1-\lambda} \varepsilon^2}\right.\\
        &\left.\quad + \frac{\kappa^2 L^2 \norm{\mx^0 - \ubmx^0}^2}{n(1-\lambda)  \varepsilon^2}  + \frac{ \sumn\norm{\nabla_x f_i(x_i^0,y_i^0)}^2}{n \varepsilon^2}}.
    \end{align*}
    Particularly, if all the agents initialize at the same starting point, i.e., $(x_i^0, y_i^0) = (x^0, y^0)$, $\forall i\in\cN$, and the underlying communication network is a complete graph, i.e., $1-\lambda^2 = 1$, the iteration complexity of DGTA becomes 
    \begin{align}
        \label{eq:dgta_complete}
        \order{\frac{\Delta_\Phi\kappa^2L}{\varepsilon^2} + \frac{\kappa D^2L^2}{\varepsilon^2}+ \frac{ \sumn\norm{\nabla_x f_i(x^0,y^0)}^2}{n \varepsilon^2}},
    \end{align}
which is similar to that of the  centralized Gradient Descent Ascent (GDA) method, i.e., $\orderi{\kappa^2 L\Delta_\Phi/\varepsilon^2 + \kappa L^2D^2/\varepsilon^2}$ \cite{lin2020gradient}.
\end{remark}




The convergence rate of DSGTA depends on the mini-batch size $b$ according to Theorem \ref{thm:dsgda}. We consider below two different choices of $b$. First, Corollary \ref{cor:comp_small_eta} states the convergence result of DSGTA under a constant batch size ($b = 1$).

\begin{corollary}[DSGTA with $b = 1$]
	\label{cor:comp_small_eta}
	Suppose Assumptions \ref{a.network}-\ref{as:sgrad} hold. Let 
	\begin{align*}
        \eta_y&=\frac{1}{\frac{9(L+\mu)}{\sqrt{1-\lambda^2}} + \sqrt{\frac{3T\sigma^2}{D^2}}} ,\; b = 1,\;
		\eta_x = \frac{1}{\sqrt{\frac{L\kappa\sigma^2 T}{n\Delta_\Phi}} + \gamma},\\
		\gamma&= \frac{120\sqrt{3}\kappa L}{(1-\lambda^2)^2} + \frac{72\sqrt{6}\kappa^2(L+\mu)}{\sqrt{1-\lambda^2}} + \frac{110\kappa^{5/4} L}{(1-\lambda^2)^{9/8}},
    \end{align*}
    for some given $T\geq 1$.
	Then, for DSGTA (Algorithm \ref{DSGTA}), we have
    \begin{align*}
        &\frac{1}{T}\sum_{t=0}^{T-1}\mE\brk{\norm{\nabla \Phi(\bar{x}^t)}^2}\leq \order{\sqrt{\frac{\Delta_\Phi L\kappa\sigma^2 }{n T }}+ \frac{ \Delta_\Phi\kappa^3 L}{(1-\lambda) T} + \frac{D^2\kappa L^2}{\sqrt{1-\lambda}T} + \frac{D\sigma L\kappa}{\sqrt{T}}\right.\\
		&\left.\quad+ \frac{\kappa^2 L^2 \norm{\mx^0 - \ubmx^0}^2}{n(1-\lambda) T}  + \frac{\sumn\norm{\nabla_x f_i(x_i^0,y_i^0)}^2}{nT} + \frac{n\Delta_\Phi\kappa L}{(1-\lambda)^3T} },
    \end{align*}
    where $\orderi{\cdot}$ hides the numerical constants.

\end{corollary}
\begin{proof}
    We have $\eta_y\mu \leq \sqrt{1-\lambda^2}/[9(\kappa+1)]$ and 
	\begin{equation}
		\label{eq:eta_inv}
		\begin{aligned}
			\frac{D^2L^2}{\eta_y\mu T} &= \frac{9D^2 \kappa L(L+\mu)}{\sqrt{1-\lambda^2} T} + \sqrt{\frac{3D^2\sigma^2 L^2\kappa^2}{T}}, \; \eta_y\kappa L\sigma^2\leq \sqrt{\frac{D^2\kappa^2 L^2\sigma^2}{3T}}\\
			\frac{\Delta_\Phi}{\eta_x T} &=  \sqrt{\frac{\Delta_\Phi L\kappa\sigma^2 }{n T }} + \frac{120\Delta_\Phi\sqrt{3}\kappa L}{(1-\lambda^2)^2 T} + \frac{72\Delta_\Phi\sqrt{6}\kappa^2(L+\mu)}{\sqrt{1-\lambda^2} T} + \frac{110\Delta_\Phi\kappa^{5/4} L}{(1-\lambda^2)^{9/8}T},\\
			\eta_x&\leq \sqrt{\frac{n\Delta_\Phi}{L\kappa\sigma^2 T}},\; \eta_x^2\leq \min\crk{\frac{n\Delta_\Phi}{L\kappa\sigma^2 T}, \frac{(1-\lambda^2)^4}{14400\kappa^2 L^2}},\;\frac{\eta_x^2\kappa^2 L^2}{\eta_y\mu}\leq \frac{(1-\lambda^2)^2}{2880\sqrt{2}}.
		\end{aligned}
	\end{equation}

	Substituting \eqref{eq:eta_inv} into \eqref{eq:dsgda_avgE} leads to 
	\begin{align*}
		&\frac{1}{T}\sum_{t=0}^{T-1}\mE\brk{\norm{\nabla \Phi(\bar{x}^t)}^2}\leq 4\sqrt{\frac{\Delta_\Phi L\kappa\sigma^2 }{n T }} + \frac{240\Delta_\Phi\sqrt{3}\kappa L}{(1-\lambda^2)^2 T} + \frac{144\Delta_\Phi\sqrt{6}\kappa^2(L+\mu)}{\sqrt{1-\lambda^2} T} + \frac{220\Delta_\Phi\kappa^{5/4} L}{(1-\lambda^2)^{9/8}T}\\
		&+ \frac{600\kappa^2 L^2 \norm{\mx^0 - \ubmx^0}^2}{n(1-\lambda^2) T} + \frac{144D^2\kappa L(L+\mu)}{\sqrt{1-\lambda^2}T} + \frac{32\sqrt{3}D\sigma L\kappa}{\sqrt{T}} + \frac{4800n\Delta_\Phi\kappa L }{(1-\lambda^2)^3T^2} + \frac{\sumn\norm{\nabla_x f_i(x_i^0,y_i^0)}^2}{nT}\\
		&\quad + \brk{\frac{162\kappa^2}{n(1-\lambda^2)} + \frac{225}{(1-\lambda^2)^3} + \frac{1}{n} }\frac{64n\Delta_\Phi\kappa L}{T},
	\end{align*}
 which yields the desired result.

\end{proof}

\begin{remark}
    Corollary \ref{cor:comp_small_eta} implies that the iteration/sample complexity of DSGTA with mini-batch size $b=1$ behaves as 
    \begin{equation}
		\label{eq:comp_iter}
		\begin{aligned}
		        &\order{\frac{\Delta_\Phi L\kappa\sigma^2 }{n \varepsilon^4 }+ \frac{D^2\sigma^2L^2\kappa^2}{\varepsilon^4}  + \frac{\kappa^2 L^2 \norm{\mx^0 - \ubmx^0}^2}{n(1-\lambda)  \varepsilon^2} + \frac{D^2\kappa L^2}{\sqrt{1-\lambda} \varepsilon^2}\right.\\
		        &\left.\quad + \frac{\sumn\norm{\nabla_x f_i(x_i^0,y_i^0)}^2}{ n\varepsilon^2}
				+ \frac{ \Delta_\Phi\kappa^3 L}{(1-\lambda) \varepsilon^2 } + \frac{n\Delta_\Phi\kappa L}{(1-\lambda)^3\varepsilon^2}}.
		    \end{aligned}
	\end{equation}
 When $\varepsilon$ is small enough, the above complexity result becomes $\order{\frac{\Delta_\Phi L\kappa\sigma^2 }{n \varepsilon^4 }+ \frac{D^2\sigma^2L^2\kappa^2}{\varepsilon^4}}$, which is independent of the graph topology. Such a result outperforms that of SGDA under a large batch size in terms of the condition number $\kappa$ \cite{lin2020gradient}.
\end{remark}



When the communication between the agents is expensive, it is often preferable to consider a large batch size $b$ in order to reduce the number of iterations.
Corollary \ref{cor:comp_large_b} states the following convergence rate of DSGTA given an arbitrary $b$.

\begin{corollary}[DSGTA with $b > 1$]
	\label{cor:comp_large_b}
	Suppose Assumptions \ref{a.network}-\ref{as:sgrad} hold. Let 
	\begin{align*}
        \eta_y&=\frac{\sqrt{1-\lambda^2}}{9(L+\mu)},\;
		\eta_x = \frac{1}{\gamma},\\
		\gamma&= \frac{120\sqrt{3}\kappa L}{(1-\lambda^2)^2} + \frac{72\sqrt{6}\kappa^2(L+\mu)}{\sqrt{1-\lambda^2}} + \frac{110\kappa^{5/4} L}{(1-\lambda^2)^{9/8}},
    \end{align*}
    for some given accuracy $\varepsilon > 0$.
	Then, for DSGTA we have
	\begin{equation}
		\label{eq:conv_large_b}
		\begin{aligned}
	        \frac{1}{T}\sum_{t=0}^{T-1}\mE\brk{\norm{\nabla \Phi(\bar{x}^t)}^2}&= \order{\frac{\Delta_\Phi\kappa L}{(1-\lambda)^2 T} + \frac{\Delta_\Phi\kappa^2L}{\sqrt{1-\lambda} T} + \frac{\kappa^2 L^2 \norm{\mx^0 - \ubmx^0}^2}{n(1-\lambda) T}\right.\\
	        &\left.\quad + \frac{D^2\kappa L^2}{\sqrt{1-\lambda}T}+ \frac{\sumn\norm{\nabla_x f_i(x_i^0,y_i^0)}^2}{nT}+ \frac{\kappa\sigma^2}{b}},
	    \end{aligned}
	\end{equation}
	where $\orderi{\cdot}$ hides numerical constants.
\end{corollary}

\begin{proof}
    We have 
	\begin{equation}
		\label{eq:eta_inv_large_b}
		\begin{aligned}
			\eta_y\mu &= \frac{\sqrt{1-\lambda^2}}{9(\kappa + 1)},\; \eta_y L =\frac{\sqrt{1-\lambda^2}\kappa}{9(\kappa + 1)}\\
			\frac{\Delta_\Phi}{\eta_x T} &= \frac{120\Delta_\Phi\sqrt{3}\kappa L}{(1-\lambda^2)^2 T} + \frac{72\Delta_\Phi\sqrt{6}\kappa^2(L+\mu)}{\sqrt{1-\lambda^2} T} + \frac{110\Delta_\Phi\kappa^{5/4} L}{(1-\lambda^2)^{9/8}T},\\
			\eta_x&\leq \frac{\sqrt{1-\lambda^2}}{72\sqrt{6}\kappa^2 (L+\mu)},\; \eta_x^2\leq \frac{(1-\lambda^2)^4}{14400\kappa^2 L^2},\;\frac{\eta_x^2\kappa^2 L^2}{\eta_y\mu}\leq \frac{(1-\lambda^2)^2}{2880\sqrt{2}},\; \frac{\eta_x^2\kappa^2 L^2}{\eta_y^2\mu^2}\leq \frac{1}{384}.
		\end{aligned}
	\end{equation}

    Substituting \eqref{eq:eta_inv_large_b} into \eqref{eq:dsgda_avgE} leads to the desired convergence rate \eqref{eq:conv_large_b}.
    %

\end{proof}

\begin{remark}
According to Corollary \ref{cor:comp_large_b}, under a batch size $b = \frac{\kappa\sigma^2}{\varepsilon^2}$ that depends on the desired accuracy level $\varepsilon$, the iteration complexity for DSGTA to obtain an $\varepsilon$-stationary point behaves as 
	\begin{align*}
		&\order{\frac{\Delta_\Phi\kappa L}{(1-\lambda)^2 \varepsilon^2} + \frac{\Delta_\Phi\kappa^2 L}{\sqrt{1-\lambda} \varepsilon^2} + \frac{\kappa^2 L^2 \norm{\mx^0 - \ubmx^0}^2}{n(1-\lambda) \varepsilon^2}+ \frac{D^2\kappa L^2}{\sqrt{1-\lambda}\varepsilon^2}+ \frac{\sumn\norm{\nabla_x f_i(x_i^0,y_i^0)}^2}{n\varepsilon^2}}.
	\end{align*}

	Therefore, the sample complexity for DSGTA is given by
	\begin{equation}
        \label{eq:comp_large_b}
	    \begin{aligned}
		&\order{\frac{\Delta_\Phi\kappa^2 L\sigma^2}{(1-\lambda)^2 \varepsilon^4} + \frac{\Delta_\Phi\kappa^3\sigma^2 L}{\sqrt{1-\lambda} \varepsilon^4} + \frac{\kappa^3\sigma^2 L^2 \norm{\mx^0 - \ubmx^0}^2}{n(1-\lambda) \varepsilon^4}+ \frac{D^2\kappa^2 L^2\sigma^2}{\sqrt{1-\lambda}\varepsilon^4}+ \frac{\sumn\norm{\nabla_x f_i(x_i^0,y_i^0)}^2\kappa\sigma^2}{n\varepsilon^4}}.
	\end{aligned}
	\end{equation}


	If all the agents initialize at the same starting point, i.e., $(x_i^0, y_i^0) = (x^0, y^0)$, $\forall i\in\cN$, and the underlying communication network is a complete graph, 
	the sample complexity for DSGTA to obtain an $\varepsilon$-stationary point behaves as 
	\begin{align*}
		&\order{\frac{\Delta_\Phi\kappa^3\sigma^2 L}{\varepsilon^4} + \frac{D^2\kappa^2 L^2\sigma^2}{\varepsilon^4}+ \frac{\sumn\norm{\nabla_x f_i(x^0,y^0)}^2\kappa\sigma^2}{n\varepsilon^4}},
	\end{align*}
which matches that of SGDA \cite{lin2020gradient} if we further assume $\sumn\normi{\nabla_x f_i(x^0,y^0)}^2\sim\orderi{n}$. 
\end{remark}

\begin{remark}
    Comparing the sample complexities \eqref{eq:comp_iter} and \eqref{eq:comp_large_b} indicates that choosing $b = 1$ can better mitigate the influence of the network. Particularly, the sample complexity in \eqref{eq:comp_iter} is dominated by $\orderi{{\Delta_\Phi L\kappa\sigma^2 }/\prti{n \varepsilon^4 }+ {D^2\sigma^2L^2\kappa^2}/\prti{\varepsilon^4}}$ for small enough $\varepsilon$. 
    This is partly due to the smaller stepsizes depending on $T$.
    Notably, such a complexity does not depend on the coefficient $1/(1-\lambda)$, which can behave as $\orderi{n^2}$ for sparse networks such as rings. Moreover, the sample complexity when $b=1$ outperforms that in \eqref{eq:comp_large_b} in terms of $\kappa$ and $1/(1-\lambda)$ in the $\orderi{\varepsilon^{-4}}$ term.

    {Nevertheless, when considering the iteration complexity, a larger batch size works better. Therefore, a tradeoff between the sampling and communication costs should be considered when choosing the best batch size $b$.}
\end{remark}

\section{Numerical Experiment}
\label{sec:sims}

This section validates the theoretical findings through  training an empirical Wasserstein robustness model \cite{sinha2017certifying} over a collection of data samples $\crk{\xi_{ij}}$ with $\ell_2$-norm attack and a penalty parameter $\gamma >0$, stated below:
\begin{equation}
\label{eq:exp}
	\begin{aligned}
		\min_{x\in\R^d} \max_{\crk{y_{ij}}_{j=1}^{m_i}\in\cY_i,i\in\cN }  \frac{1}{n}\sum_{i = 1}^{n} \crk{   \frac{1}{m_i}\sum_{j=1}^{m_i}\bigg[\ell(x; y_{ij}) - \gamma\|y_{ij} - \xi_{ij}\|^2_2\bigg]  },
	\end{aligned}
\end{equation}
where $\xi_{ij}$ is the $j$-the samples of agent $i$ and {$\ell:\R^d\times \cY_i\rightarrow\R$ is a neural network classifier.} In particular, {the classifier comprises three convolutional layers with filter sizes $8\times 8$, $6\times 6$ and $5\times 5$, each with a ELU activation, followed by a fully connected layer and softmax output.}
As demonstrated in \cite{sinha2017certifying}, the component function $\brki{\ell(x; y_{ij}) - \gamma\normi{y_{ij} - \xi_{ij}}^2_2}$ is strongly concave in $y_{ij}$ with sufficiently large $\gamma >0$ for any $i\in\cN$ (e.g., $\gamma = 0.4$ or $\gamma = 1.3$ for small and large adversarial perturbations, respectively). In this experiment, we set $\gamma = 0.4$ to ensure NC-SC objective functions.



\begin{figure}[htbp]
    \centering
    \subfloat[Comparison between DGTA and GT/DA with various $\eta_x$.]{\includegraphics[width=0.48\textwidth]{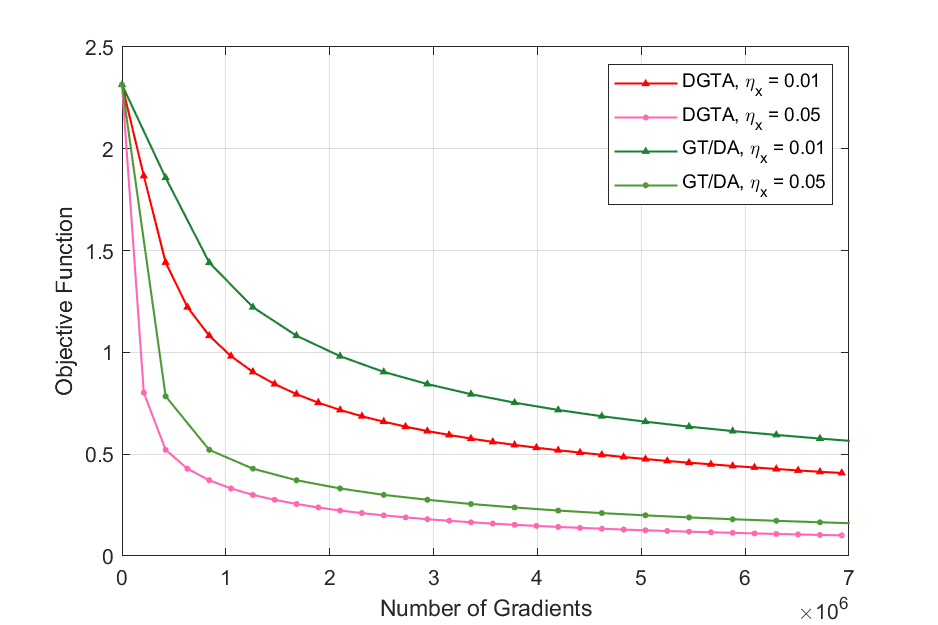}\label{fig:det}}
    \subfloat[DSGTA with $b=8$ and various $\eta_x$.]{\includegraphics[width=0.48\textwidth]{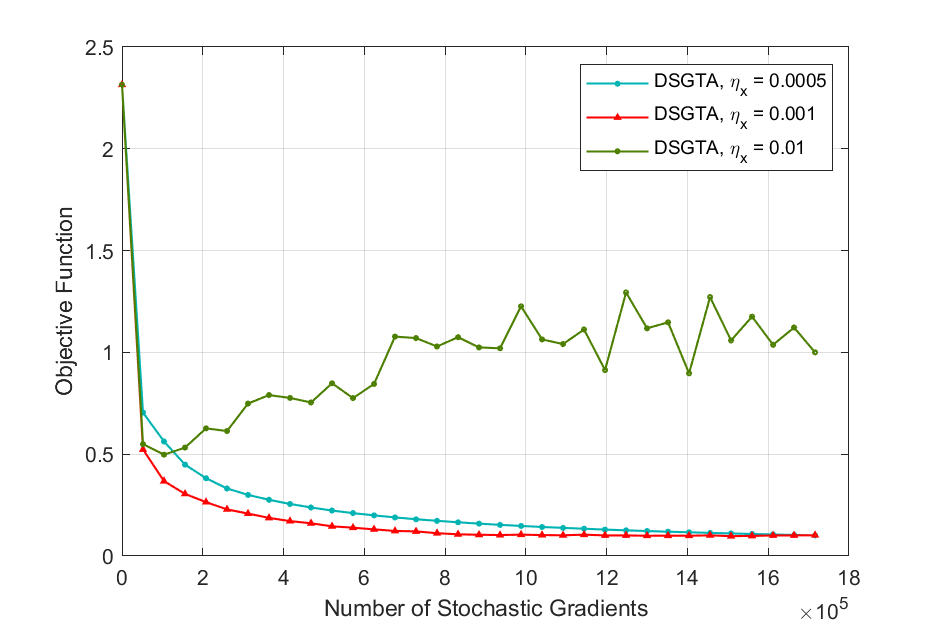}\label{fig:sto}}
    \caption{Performance comparison of DSGTA, DGTA, and GT/DA under the same $\eta_y=0.01$ and varying stepsizes $\eta_x$. The horizontal axis (``Number of Gradients''  and ``Number of Stochastic Gradients'') shows the iteration progress, where each gradient  or stochastic gradient evaluation with respect to either variable $x$ or $y$ increments the count by one.}
    \label{fig:graph}
\end{figure}

We conduct the experiments on the MNIST dataset \cite{lecun2010mnist} and  consider both the full gradient and the stochastic gradient settings over a ring network of 24 agents. The reported results are averaged over three independent runs. We use the function value $\frac{1}{n}\sum_{i = 1}^{n} \crki{ \frac{1}{m_i}\sum_{j=1}^{m_i}[\ell(\bar{x}^t; y_{ij}^t) - \gamma\|y_{ij}^t - \xi_{ij}\|^2_2]}$ 
as a metric to demonstrate the algorithms' performance. 
In particular, we  compare DGTA and GT/DA \cite{9054056}  under the full gradient setting in Figure~\ref{fig:det} with varying stepsizes $\eta_x\in\crki{0.01, 0.05}$ and fixed $\eta_y = 0.01$. The GT/DA  method incorporates three inner-loop updates for the variables $y_i$. Figure~\ref{fig:det} indicates that DGTA outperforms GT/DA under the same stepsize.

 Under the stochastic gradient setting (Figure~\ref{fig:sto}), we illustrate the performance of DSGTA under stepsizes $\eta_x\in\crki{0.0005,0.001, 0.01}$ and fixed $\eta_y = 0.01$. Figure~\ref{fig:sto} demonstrates that proper choice of $\eta_x (\ll\eta_y)$ is critical for the success of the DSGTA method in practice. It is worth noting that by employing stochastic gradients, the convergence with respect to the number of samples can be greatly accelerated, since the full gradient computation requires.

\section{Conclusion}
\label{sec:conclusion}

This paper focuses on addressing the distributed min-max optimization problem over networked agents under a novel formulation which enables wide applications. 
The proposed algorithms, Distributed Gradient Tracking Ascent (DGTA) and Distributed Stochastic Gradient Tracking Ascent
(DSGTA), leverage gradient tracking technique to achieve the comparable theoretical guarantees as their centralized counterparts (up to constant factors related to the communication graph). In particular, We demonstrate that DGTA achieves an iteration complexity of $\orderi{\kappa^2\varepsilon^{-2}}$, and DSGTA enjoys a sample complexity of $\orderi{\kappa^3\varepsilon^{-4}}$ for nonconvex strongly concave objective functions. 
Experimental results also corroborate to such theoretical findings. 

\bibliographystyle{siamplain}
\bibliography{references_all}

\newpage
\appendix

\section*{Appendix}

\section{Proofs}

\subsection{Proof of Lemma \ref{lem:dsgda_descent}}
\label{app:dsgda_descent}

To facilitate the analysis, we define $\cF_t, t= 0,1,\ldots$ to be the $\sigma$-algebra generated by $\crk{\bxi^0_i, \bxi^1_i,\cdots, \bxi^{t-1}_i}_{i=1}^n$ according to Line \ref{line:dsgta_xi} in Algorithm \ref{DSGTA}. 
	Note that $\bar{x}^{t + 1} = \bar{x}^t - \eta_x \bg^t$ due to Line~\ref{line:dsgta_x} in Algorithm~\ref{DSGTA} and Assumption \ref{a.network}. We then apply the descent lemma to $\Phi(\cdot)$ and obtain
	\begin{equation}
		\label{eq:dsgda_descent_s1}
		\begin{aligned}
			&\condE{\Phi(\bar{x}^{t + 1})}{\cF_{t}} \leq \Phi(\bar{x}^t) - \eta_x \inpro{\nabla \Phi(\bar{x}^t), \frac{1}{n}\sumn\nabla f_i(x_i^t,y_i^t)} + \kappa L\eta_x^2\condE{\norm{\bg^t}^2}{\cF_{t}}\\
			&=  \Phi(\bar{x}^t) -\frac{\eta_x}{2}\prt{-\norm{\nabla\Phi(\bar{x}^t) - \frac{1}{n}\sumn\nabla f_i(x_i^t,y_i^t)}^2 + \norm{\frac{1}{n}\sumn\nabla f_i(x_i^t,y_i^t)}^2 + \norm{\nabla\Phi(\bar{x}^t)}^2}\\
			&\quad + \kappa L\eta_x^2\crk{\condE{\norm{\bg^t - \frac{1}{n}\sumn\nabla f_i(x_i^t,y_i^t)}^2}{\cF_{t}} + \norm{\frac{1}{n}\sumn\nabla f_i(x_i^t,y_i^t)}^2}\\
			&= \Phi(\bar{x}^t) - \frac{\eta_x}{2}\norm{\nabla \Phi(\bar{x}^t)}^2 - \frac{\eta_x}{2}\prt{1 - 2\kappa L\eta_x}\norm{\frac{1}{n}\sumn\nabla f_i(x_i^t,y_i^t)}^2 \\
			&\quad + \frac{\eta_x}{2} \norm{\nabla \Phi(\bar{x}^t) - \frac{1}{n}\sum_{i=1}^n\nabla_x f_i(x_i^t, y_i^{t})}^2 + \kappa L\eta_x^2\condE{\norm{\bg^t - \frac{1}{n}\sumn\nabla f_i(x_i^t,y_i^t)}^2}{\cF_{t}} \\
			&\leq \Phi(\bar{x}^t) - \frac{\eta_x}{2}\norm{\nabla\Phi(\bar{x}^t)}^2 -\frac{\eta_x}{2}\prt{1 - 2\kappa L\eta_x}\norm{\frac{1}{n}\sumn\nabla_x f_i(x_i^t,y_i^t)}^2 + \frac{\eta_x L^2}{n}\sum_{i=1}^n \norm{\hat{y}_i(\bar{x}^t) - y_i^t}^2\\
			&\quad + \frac{\eta_x L^2}{n}\sumn\norm{\bar{x}^t - x_i^t}^2 + \frac{\kappa L\eta_x^2\sigma^2}{nb},
		\end{aligned}
	\end{equation}
	where the last inequality holds due to Assumption~\ref{as:sgrad} and 
	\begin{equation}
		\label{eq:phi_avg}
		\begin{aligned}
			\norm{\nabla \Phi(\bar{x}^t) - \frac{1}{n}\sum_{i=1}^n\nabla_x f_i(x_i^t, y_i^{t})}^2 &= \norm{\frac{1}{n}\sumn\nabla_x f_i(\bar{x}^t, \hat{y}_i(\bar{x}^2)) - \frac{1}{n}\sum_{i=1}^n\nabla_x f_i(x_i^t, y_i^{t})}^2\\
			&\leq \frac{2L^2}{n}\sumn\brk{\norm{\bar{x}^t - x_i^t}^2 + \norm{\hat{y}_i(\bar{x}^t) - y_i^t}^2}.
		\end{aligned}
	\end{equation}

	Taking the full expectation on both sides of \eqref{eq:dsgda_descent_s1} yields the desired result.


\subsection{Proof of Lemma \ref{lem:dsgda_delta}}
\label{app:dsgda_delta}

	Denote $\delta_i^t := \normi{\hat{y}_i(\bar{x}^t) - y_i^t}^2$. Recall that $f_i(\bar{x}^t, \cdot)$ is $\mu$-strongly concave. Then, 
	\begin{equation}
		\label{eq:dsgda_delta_s1}
		\begin{aligned}
			f_i(\bar{x}^t, \hat{y}_i(\bar{x}^t))
			&\leq f_i(\bar{x}^t, y_i^t) + \inpro{\nabla_y f_i(\bar{x}^t, y_i^t), \hat{y}_i(\bar{x}^t) - y_i^t} - \frac{\mu}{2}\delta_i^t.
		\end{aligned}
	\end{equation}

	Denote $z_i^t := y_i^t + \theta \uh_{y,i}^t$, where $\uh_{y,i}^t = \sum_{j=1}^b h_{y,i}(x_i^t, y_i^t,\xi_{i,j}^t) / b$ is the noisy gradient obtained by sampling $\crki{\xi_{i,j}^t}_{j=1}^b$. From Remark \ref{rem:smooth} and the descent lemma, we have  
	\begin{equation}
		\label{eq:dsgda_delta_s2}
		\begin{aligned}
			\condE{f_i(\bar{x}^t,y_i^t)}{\cF_t}
			&\leq \condE{f_i(\bar{x}^t, z_i^t)}{\cF_t}  - \theta \condE{\inpro{ \nabla_y f_i(\bar{x}^t,z_i^t), \uh_{y,i}^t}}{\cF_t}\\
			&\quad +  L\theta^2\condE{\norm{\uh_{y,i}^t}^2}{\cF_t}.
		\end{aligned}
	\end{equation}

	Note that $\hat{y}_i(\cdot) = \arg\max_{y\in\mY_i} f_i(\cdot, y)$. We have $f_i(\bar{x}^t, z_i^t)\leq f_i(\bar{x}^t, \hat{y}_i(\bar{x}^t))$. Substituting such a relation into \eqref{eq:dsgda_delta_s2} and recalling \eqref{eq:dsgda_delta_s1} lead to 

	\begin{equation}
		\label{eq:dsgda_delta_s3}
		\begin{aligned}
			\condE{f_i(\bar{x}^t, y_i^t)}{\cF_t}&\leq \condE{f_i(\bar{x}^t, y_i^t)}{\cF_t} + \condE{\inpro{\nabla_y f_i(\bar{x}^t, y_i^t), \hat{y}_i(\bar{x}^t) - y_i^t}}{\cF_t} - \frac{\mu}{2}\condE{\delta_i^t}{\cF_t}\\
			&\quad - \theta \condE{\inpro{ \nabla_y f_i(\bar{x}^t,z_i^t), \uh_{y,i}^t}}{\cF_t}+  L\theta^2\condE{\norm{ \uh_{y,i}^t}^2}{\cF_t}.
		\end{aligned}
	\end{equation}

	Rearranging the terms in \eqref{eq:dsgda_delta_s3} and invoking Assumption \ref{as:sgrad} yield
	\begin{equation}
		\label{eq:dsgda_delta_up}
		\begin{aligned}
			&0 \leq \condE{\inpro{\nabla_y f_i(\bar{x}^t, y_i^t), \hat{y}_i(\bar{x}^t) - y_i^t}}{\cF_t} - \theta \condE{\inpro{ \nabla_y f_i(\bar{x}^t,z_i^t), \uh_{y,i}^t}}{\cF_t}\\
			&\quad +  \frac{L\theta^2\sigma^2}{b} + L\theta^2\norm{\nabla_y f_i(x_i^t,y_i^t)}^2 -\frac{\mu}{2}\condE{\delta_i^t}{\cF_t}\\
			&= \condE{\inpro{\nabla_y f_i(\bar{x}^t, y_i^t), \hat{y}_i(\bar{x}^t) - z_i^t}}{\cF_t} + \theta\condE{\inpro{\nabla_y f_i(\bar{x}^t, y_i^t), \uh_{y,i}^t}}{\cF_t} -\frac{\mu}{2}\condE{\delta_i^t}{\cF_t}\\
			&\quad - \theta \condE{\inpro{ \nabla_y f_i(\bar{x}^t,z_i^t), \uh_{y,i}^t}}{\cF_t} +  \frac{L\theta^2\sigma^2}{b} + L\theta^2\norm{\nabla_y f_i(x_i^t,y_i^t)}^2\\
			&= \condE{\inpro{\nabla_y f_i(\bar{x}^t, y_i^t) - \nabla_y f_i(x_i^t, y_i^t), \hat{y}_i(\bar{x}^t) - z_i^t}}{\cF_t}  + \condE{\inpro{\nabla_y f_i(x_i^t, y_i^t), \hat{y}_i(\bar{x}^t) - y_i^t - \theta \uh_{y,i}^t}}{\cF_t}\\
			&\quad + \theta \condE{\inpro{ \nabla_y f_i(\bar{x}^t, y_i^t) - \nabla_y f_i(\bar{x}^t,z_i^t), \uh_{y,i}^t}}{\cF_t} +  \frac{L\theta^2\sigma^2}{b} + L\theta^2\norm{\nabla_y f_i(x_i^t,y_i^t)}^2-\frac{\mu}{2}\condE{\delta_i^t}{\cF_t}\\
			&= \condE{\inpro{\nabla_y f_i(\bar{x}^t, y_i^t) - \nabla_y f_i(x_i^t, y_i^t), \hat{y}_i(\bar{x}^t) - z_i^t}}{\cF_t}  + \condE{\inpro{\nabla_y f_i(x_i^t, y_i^t), \hat{y}_i(\bar{x}^t) - y_i^t}}{\cF_t}-\frac{\mu}{2}\condE{\delta_i^t}{\cF_t}\\
			&\quad + \theta \condE{\inpro{\nabla_y f_i(\bar{x}^t, y_i^t)- \nabla_y f_i(\bar{x}^t,z_i^t), \uh_{y,i}^t}}{\cF_t} +  \frac{L\theta^2\sigma^2}{b} + \prt{L\theta^2-\theta}\norm{\nabla_y f_i(x_i^t,y_i^t)}^2,
		\end{aligned}
	\end{equation}
	where we invoked Assumption~\ref{as:sgrad} that $\condEi{\uh_{y,i}^t}{\cF_t} = \nabla_y f_i(x_i^t,y_i^t)$ in the last equality.

	Then it remains to estimate the inner products in \eqref{eq:dsgda_delta_up}. We start by bounding the last inner product in \eqref{eq:dsgda_delta_up}. Recall that $\nabla_y f_i(\bar{x}^t,\cdot)$ is $2L$-Lipschitz continuous. By the Cauchy-Schwarz inequality, we have 
	\begin{equation}
		\label{eq:dsgda_delta_inner0}
		\begin{aligned}
			&\theta \condE{\inpro{\nabla_y f_i(\bar{x}^t, y_i^t)- \nabla_y f_i(\bar{x}^t,z_i^t), \uh_{y,i}^t}}{\cF_t}\\
			&\leq \theta \condE{\norm{\nabla_y f_i(\bar{x}^t, y_i^t)- \nabla_y f_i(\bar{x}^t,z_i^t)}\norm{\uh_{y,i}^t }}{\cF_t}\\
			&\leq 2\theta^2 L \condE{\norm{\uh_{y,i}^t}^2}{\cF_t}\\
			&\leq 2\theta^2 L\prt{ \frac{\sigma^2}{b} + \norm{\nabla_y f_i(x_i^t, y_i^t)}^2}.
		\end{aligned}
	\end{equation}

	For the first inner product in \eqref{eq:dsgda_delta_up}, we have from Young's inequality that for any $q>0$,
	\begin{equation}
		\label{eq:dsgda_delta_inner1}
		\begin{aligned}
			& \condE{\inpro{\nabla_y f_i(\bar{x}^t, y_i^t) - \nabla_y f_i(x_i^t, y_i^t), \hat{y}_i(\bar{x}^t) - z_i^t}}{\cF_t}\\
			&\leq \frac{q}{2}\norm{\nabla_y f_i(\bar{x}^t, y_i^t) - \nabla_y f_i(x_i^t, y_i^t)}^2 + \frac{1}{2q}\condE{\norm{\hat{y}_i(\bar{x}^t) - y_i^t - \theta \uh_{y,i}^t}^2}{\cF_t}\\
			&\leq qL^2\norm{\bar{x}^t - x_i^t}^2 + \frac{1}{q}\condE{\delta_i^t}{\cF_t} + \frac{\theta^2}{q}\prt{\norm{\nabla_y f_i(x_i^t, y_i^t)}^2 + \frac{\sigma^2}{b}}.
		\end{aligned}
	\end{equation}
	
	Note that $y_i^{t + 1} = \proj_{\mY_i}\prt{z_i^t}$ when choosing $\theta = \eta_y$ in $z_i^t = y_i^t + \theta h_{y,i}^t$. For the remaining inner product in \eqref{eq:dsgda_delta_up}, we have
	\begin{equation}
		\label{eq:dsgda_delta_inner2}
		\begin{aligned}
			& 2\eta_y \condE{\inpro{\nabla_y f_i(x_i^t, y_i^t), \hat{y}_i(\bar{x}^t) - y_i^t}}{\cF_t} = 2\eta_y\condE{\inpro{ \uh_{y,i}^t, \hat{y}_i(\bar{x}^t) - y_i^t}}{\cF_t}\\
			&\quad + 2\eta_y\condE{\inpro{\nabla_y f_i(x_i^t, y_i^t) - \uh_{y,i}^t, \hat{y}_i(\bar{x}^t) - y_i^t}}{\cF_t}\\
			&= 2\condE{\inpro{z_i^t - y_i^t, \hat{y}_i(\bar{x}^t) - y_i^t}}{\cF_t}\\
			&= \condE{\norm{z_i^t - y_i^t}^2}{\cF_t} + \norm{\hat{y}_i(\bar{x}^t) - y_i^t}^2 - \condE{\norm{\hat{y}_i(\bar{x}^t) - z_i^t}^2}{\cF_t}.
		\end{aligned}
	\end{equation}

	The last term in \eqref{eq:dsgda_delta_inner2} can be bounded as follows:
	\begin{equation}
		\label{eq:dsgda_delta_inner2_s1}
		\begin{aligned}
			\condE{\norm{\hat{y}_i(\bar{x}^t) - z_i^t}^2}{\cF_t} &= \condE{\norm{\hat{y}_i(\bar{x}^t) - y_i^{t + 1}}^2}{\cF_t} + \condE{\norm{y_i^{t + 1} - z_i^t}^2}{\cF_t}\\
			&\quad + 2\condE{\inpro{\hat{y}_i(\bar{x}^t) - y_i^{t + 1}, y_i^{t + 1} - z_i^t}}{\cF_t}\\
			&\geq \condE{\norm{\hat{y}_i(\bar{x}^t) - y_i^{t + 1}}^2}{\cF_t} + \condE{\norm{y_i^{t + 1} - z_i^t}^2}{\cF_t},
		\end{aligned}
	\end{equation}
	where we invoke the second projection Theorem \cite{beck2017first}, i.e., given a closed convex set $\mY$,
	\begin{align}
		\label{eq:second_proj2}
		\inpro{y - \proj_{\mY}(u), \proj_{\mY}(u) - u}\geq 0, \ \forall y\in\mY.
	\end{align}
	The inequality \eqref{eq:dsgda_delta_inner2_s1} holds by choosing $y = \hat{y}_i(\bar{x}^t)$, $u = z_i^t$, and $\mY = \mY_i$ in \eqref{eq:second_proj2}. Therefore, \eqref{eq:dsgda_delta_inner2} becomes 
	\begin{equation}
		\label{eq:dsgda_delta_inner2_s2}
		\begin{aligned}
			&\condE{\inpro{\nabla_y f_i(x_i^t, y_i^t), \hat{y}_i(\bar{x}^t) - y_i^t}}{\cF_t}\\
			&\leq \frac{1}{2\eta_y}\crk{\eta_y^2\prt{\frac{\sigma^2}{b} + \norm{\nabla_y f_i(x_i^t, y_i^t)}^2} + \condE{\delta_i^t}{\cF_t} - \condE{\norm{y_i^{t + 1} - \hat{y}_i(\bar{x}^t)}^2}{\cF_t}}.
		\end{aligned}
	\end{equation}

	Substituting \eqref{eq:dsgda_delta_inner0}, \eqref{eq:dsgda_delta_inner1}, and \eqref{eq:dsgda_delta_inner2_s2} into \eqref{eq:dsgda_delta_up} and letting $\theta = \eta_y$ lead to 
	\begin{equation}
		\label{eq:dsgda_delta_s4}
		\begin{aligned}
			&\condE{\norm{y_i^{t + 1} - \hat{y}_i(\bar{x}^t)}^2}{\cF_t} \leq \prt{1 - \eta_y\mu + \frac{2\eta_y}{q}}\condE{\delta_i^t}{\cF_t} + 2qL^2\eta_y \norm{\bar{x}^t - x_i^t}^2\\
			&\quad -\eta_y^2 \prt{2 - 2\eta_y L - 1-4\eta_y L - \frac{2\eta_y}{q}}\norm{\nabla_y f_i(x_i^t, y_i^t)}^2 + \prt{1 + \frac{2\eta_y}{q} + 4\eta_y L + 2\eta_y L}\frac{\eta_y^2\sigma^2}{b}\\
			&\leq \prt{1 - \frac{\eta_y\mu}{2}}\delta_i^t + \frac{8 L^2\eta_y}{\mu}\norm{\bar{x}^t - x_i^t}^2 + \frac{2\eta_y^2\sigma^2}{b},
		\end{aligned}
	\end{equation}
	where the last inequality holds by letting $q = 4/\mu$ and $\eta_y\leq 1/[8(\mu + L)]$. 
 
 We are now ready to derive the recursion for $\delta_i^t$:
	\begin{align}
		&\condE{\delta_i^{t+1}}{\cF_t} = \condE{\norm{\hat{y}_i(\bar{x}^{t + 1}) - y_i^{t + 1}}^2}{\cF_t}\nonumber\\
		&\leq \prt{1 + q}\condE{\norm{y_i^{t + 1} - \hat{y}_i(\bar{x}^t)}^2}{\cF_t} + \prt{1 + q^{-1}}\condE{\norm{\hat{y}_i(\bar{x}^{t+1}) - \hat{y}_i(\bar{x}^t)}^2}{\cF_t}\nonumber\\
		&\leq (1 + q)\brk{\prt{1 - \frac{\eta_y\mu}{2}}\condE{\delta_i^t}{\cF_t} + \frac{8L^2\eta_y}{\mu}\norm{\bar{x}^t - x_i^t}^2 + \frac{2\eta_y^2\sigma^2}{b}} + \prt{1 + q^{-1}}\kappa^2\condE{\norm{\bar{x}^{t + 1} - \bar{x}^t}^2}{\cF_t}\label{eq:dsgda_delta_s5}\\
		&\leq \prt{1 - \frac{\eta_y\mu}{4}}\condE{\delta_i^t}{\cF_t} + \frac{9L^2\eta_y}{\mu}\norm{\bar{x}^t - x_i^t}^2 + \frac{3\eta_y^2\sigma^2}{b} + \frac{4\kappa^2\eta_x^2}{\eta_y\mu}\norm{\frac{1}{n}\sumn\nabla_x f_i(x_i^t,y_i^t)}^2 + \frac{4\kappa^2\eta_x^2\sigma^2}{\eta_y\mu nb},\label{eq:dsgda_delta_s6}
	\end{align}
	where we use the $\kappa$-Lipschitz property of $\hat{y}_i(\cdot)$ in \eqref{eq:dsgda_delta_s5} and choose $q=\mu\eta_y/(4-2\eta_y\mu)$. Summing  over all the agents and taking the full expectation on both sides of \eqref{eq:dsgda_delta_s6} lead to the desired result.

\subsection{Proof of Lemma \ref{lem:dsgda_cons0}}
\label{app:dsgda_cons0}

	Recall that $\mx^{t + 1}  = W\prti{\mx^t - \eta_x \mg^t}$ and $\normi{W - \1\1^{\T}/n}_2\leq \lambda<1$. According to Young's inequality, for any $q>0$, we have 
	\begin{equation}
		\label{eq:dsgda_cons0_s1}
		\begin{aligned}
			\condE{\norm{\mx^{t + 1} - \ubmx^{t + 1}}^2}{\cF_t} &\leq (1 + q)\lambda^2\norm{\mx^t - \ubmx^t}^2 + (1+q^{-1})\lambda^2\eta_x^2\condE{\norm{\mg^t - \ubg^t}^2}{\cF_t}\\
			&\leq  \frac{1 + \lambda^2}{2}\norm{\mx^t - \ubmx^t}^2 + \frac{2\eta_x^2}{1-\lambda^2}\condE{\norm{\mg^t - \ubg^t}^2}{\cF_t},
		\end{aligned}
	\end{equation}
	where the last inequality holds by letting $q = (1-\lambda^2)/(2\lambda^2)$. Taking the full expectation on both sides of \eqref{eq:dsgda_cons0_s1} yields the desired result.

\subsection{Proof of Lemma \ref{lem:dsgda_gtc}}
\label{app:dsgda_gtc}

	Denote $\md^t:= \h_{x}^t - \nabla_x F(\mx^t,\my^t)$ and $\cT:= \md^{t + 1} - \md^t +\nabla_x  F(\mx^{t + 1},\my^{t + 1})-\nabla_x  F(\mx^{t },\my^{t })$. From Line \ref{line:gt} in Algorithm \ref{DSGTA}, we have 
	\begin{equation}
		\label{eq:dsgda_gt_s1}
		\begin{aligned}
			&\mg^{t + 1} - \ubg^{t + 1} = \prt{W - \frac{\1\1^{\T}}{n}}\left\{\mg^t - \ubg^t + \cT \right\}.
		\end{aligned}
	\end{equation}

	Then,
	\begin{equation}
		\label{eq:dsgda_gt_norm}
		\begin{aligned}
			\condE{\norm{\mg^{t + 1} - \ubg^{t + 1}}^2}{\cF_{t + 1}} &\leq \lambda^2\norm{\mg^t - \ubg^t}^2 + \lambda^2\condE{\norm{\cT}^2}{\cF_{t + 1}}\\
			&\quad + 2\condE{\inpro{\prt{W - \frac{\1\1^{\T}}{n}}\prt{\mg^t - \ubg^t},\prt{W - \frac{\1\1^{\T}}{n}}\cT}}{\cF_{t + 1}}.
		\end{aligned}
	\end{equation}

	According to Assumption \ref{as:sgrad} and Young's inequality for any $q>0$, we have 
	\begin{equation}
		\label{eq:dsgda_gt_inner0}
		\begin{aligned}
			&2\condE{\inpro{\prt{W - \frac{\1\1^{\T}}{n}}\prt{\mg^t - \ubg^t},\prt{W - \frac{\1\1^{\T}}{n}}\cT}}{\cF_{t + 1}}\\
			&\leq -2\inpro{\prt{W - \frac{\1\1^{\T}}{n}}\prt{\mg^t - \ubg^t},\prt{W - \frac{\1\1^{\T}}{n}}\md^t}\\
			&\quad + q\lambda^2\norm{\mg^t - \ubg^t}^2 + q^{-1}\lambda^2 L^2\prt{\norm{\mx^{t + 1} - \mx^t}^2 + \condE{\norm{\my^{t + 1} - \my^t}^2}{\cF_{t + 1}}}.
		\end{aligned}
	\end{equation}

	We next consider bounding the term $\condEi{\normi{\cT}^2}{\cF_{t + 1}}$:
	\begin{equation}
		\label{eq:dsgda_gt_norm2}
		\begin{aligned}
			\condE{\norm{\cT}^2}{\cF_{t + 1}} &\leq \frac{3n\sigma^2}{b} + 3\norm{\md^t}^2 + 3L^2\prt{\norm{\mx^{t + 1} - \mx^t}^2 + \condE{\norm{\my^{t + 1} - \my^t}^2}{\cF_{t + 1}}}.
		\end{aligned}
	\end{equation}

        Similar to the procedures in \cite[Lemma 7]{pu2021distributed}, we can show that 
        \begin{equation}
            \label{eq:dsgda_gt_inner_k}
            \begin{aligned}
                \condE{-2\inpro{\prt{W - \frac{\1\1^{\T}}{n}}\prt{\mg^t - \ubg^t},\prt{W - \frac{\1\1^{\T}}{n}}\md^t}}{\cF_t}\leq \frac{\sigma^2}{b}.
            \end{aligned}
        \end{equation}
        
	Substituting \eqref{eq:dsgda_gt_inner0}-\eqref{eq:dsgda_gt_inner_k} into \eqref{eq:dsgda_gt_norm} and invoking the tower property, we have 
	\begin{equation}
		\label{eq:dsgda_gt_cFt}
		\begin{aligned}
			&\condE{\norm{\mg^{t + 1} - \ubg^{t + 1}}^2}{\cF_t} = \condE{\condE{\norm{\mg^{t + 1} - \ubg^{t + 1}}^2}{\cF_{t+1}}}{\cF_t}\\
			&\leq (1 + q)\lambda^2\condE{\norm{\mg^t - \ubg^t}^2}{\cF_t} + (3+q^{-1})L^2\prt{\condE{\norm{\mx^{t + 1} - \mx^t}^2}{\cF_t} + \condE{\norm{\my^{t + 1} - \my^t}^2}{\cF_t}} + \frac{7n\sigma^2}{b}.
		\end{aligned}
	\end{equation}

	We then bound the terms $\condEi{\normi{\mx^{t + 1} - \mx^t}^2}{\cF_t}$ and $\condEi{\normi{\my^{t + 1} - \my^t}^2}{\cF_t}$. First, we have 
	\begin{equation}
		\label{eq:dsgda_gt_xs}
		\begin{aligned}
			\frac{1}{3}\condE{\norm{\mx^{t + 1} - \mx^t}^2}{\cF_t}&\leq \condE{\norm{\mx^{t + 1} - \ubmx^{t + 1}}^2}{\cF_t} + \norm{\mx^t - \ubmx^t}^2 + n\eta_x^2\condE{\norm{\ubg^t}^2}{\cF_t}\\
			&\leq 2\norm{\mx^t - \ubmx^t}^2 + \frac{2\eta_x^2}{1-\lambda^2}\condE{\norm{\mg^t - \ubg^t}^2}{\cF_t} + n\eta_x^2\condE{\norm{\ubg}^2}{\cF_t},
		\end{aligned}
	\end{equation}
	where the last inequality holds due to \eqref{eq:dsgda_cons0_s1}.

	Recall that $\delta^t=\sumn\normi{\hat{y}_i(\bar{x}^t) - y_i^t}^2$ and $\hat{y}_i(\cdot)$ is $\kappa$-Lipschitz continuous. We have 
	\begin{equation}
		\label{eq:dsgda_gt_ys}
		\begin{aligned}
			\frac{1}{3}\condE{\norm{\my^{t + 1} - \my^t}^2}{\cF_t} &\leq \condE{\delta^{t+1}}{\cF_t} + \condE{\delta^t}{\cF_t} + \kappa^2 \eta_x^2n\condE{\norm{\ubg^t}^2}{\cF_t}\\
			&\leq 2\condE{\delta^t}{\cF_t} + \frac{9\eta_y L^2}{\mu}\norm{\mx^t - \ubmx^t}^2 + \frac{3n\eta_y^2\sigma^2}{b} + \frac{(4+\eta_y\mu)\kappa^2\eta_x^2n}{\eta_y\mu}\condE{\norm{\ubg^t}^2}{\cF_t}.
		\end{aligned}
	\end{equation}

	Substituting \eqref{eq:dsgda_gt_xs} and \eqref{eq:dsgda_gt_ys} into \eqref{eq:dsgda_gt_cFt} and letting $q = (1-\lambda^2)/(2\lambda^2)$ in \eqref{eq:dsgda_gt_cFt} lead to 
	\begin{equation}
		\label{eq:dsgda_gt_cFt2}
		\begin{aligned}
			&\condE{\norm{\mg^{t + 1} - \ubg^{t + 1}}^2}{\cF_t} \leq \prt{\frac{1+\lambda^2}{2} + \frac{18\eta_x^2 L^2}{(1-\lambda^2)^2}}\condE{\norm{\mg^t - \ubg^t}^2}{\cF_t} + \frac{18 L^2}{1-\lambda^2}\condE{\delta^t}{\cF_t}\\
			&\quad + \prt{2 + \frac{9\eta_y L^2}{\mu}}\frac{9L^2}{1-\lambda^2}\norm{\mx^t - \ubmx^t}^2 + \frac{9\eta_x^2 nL^2}{1-\lambda^2}\prt{1 + \frac{(4 + \eta_y\mu)\kappa^2}{\eta_y\mu}}\condE{\norm{\ubg^t}^2}{\cF_t} + \prt{2 + \frac{9\eta_y^2 L^2}{1-\lambda^2}}\frac{n\sigma^2}{b}.
		\end{aligned}
	\end{equation}

	Note that $\condE{\normi{\ubg^t}^2}{\cF_t}\leq \sigma^2/(nb) + \normi{\sumn\nabla_x f_i(x_i^t,y_i^t)/n}^2$.
	Letting $\eta_x \leq (1-\lambda^2)^{2/3}/(6\sqrt{2} L)$ in \eqref{eq:dsgda_gt_cFt2} and taking the full expectation lead to the desired result.

\subsection{Proof of Lemma \ref{lem:dsgda_cL}}
\label{app:dsgda_At}

	 The following procedures demonstrate how to obtain the coefficients in \eqref{eq:dsgda_lya_full}. We consider the following candidate:
	\begin{equation}
		\label{eq:dsgda_lya}
		\begin{aligned}
			\cL^t = \Phi(\bar{x}^t) - \Phi^* + \eta_x\cC_1\norm{\mx^t - \ubmx^t}^2 + \eta_x\cC_2\delta^t + \eta_x^3\cC_3\norm{\mg^t - \ubg^t}^2,
		\end{aligned}
	\end{equation}
	where $\cC_1$-$\cC_3$ are some positive constants to be determined later. 

	In light of Lemmas \ref{lem:dsgda_descent}-\ref{lem:dsgda_gtc} and \eqref{eq:dsgda_lya}, we have 
	\begin{equation}
		\label{eq:dsgda_lya1_0}
		\begin{aligned}
			\mE\cL^{t + 1} &\leq  \mE\brk{\Phi(\bar{x}^t) - \Phi^*} - \frac{\eta_x}{2}\mE\brk{\norm{\nabla \Phi(\bar{x}^t)}^2} + \frac{\eta_x^2 L\kappa\sigma^2}{nb} + \frac{3\eta_x\eta_y^2n\cC_2\sigma^2}{b}\\
			&\quad + \brk{\frac{L^2}{n} + \frac{1 + \lambda^2}{2}\cC_1 + \frac{9 L^2\eta_y\cC_2}{\mu} + \prt{2 + \frac{9\eta_y L^2}{\mu}}\frac{9\eta_x^2 L^2\cC_3}{1-\lambda^2}}\eta_x\mE\brk{\norm{\mx^t - \ubmx^t}^2}\\
			&\quad + \brk{\frac{L^2}{n} + \prt{1 - \frac{\eta_y\mu}{4}}\cC_2 + \frac{18\eta_x^2 L^2\cC_3}{1-\lambda^2} }\eta_x\mE\brk{\delta^t}\\
			&\quad + \brk{\frac{2\lambda^2\cC_1}{1-\lambda^2} + \frac{3+\lambda^2}{4}\cC_3}\eta_x^3\mE\brk{\norm{\mg^t - \ubg^t}^2}\\
			&\quad - \frac{\eta_x}{2}\brk{1 - 2\kappa L\eta_x - \frac{8n\kappa^2\eta_x^2\cC_2}{\eta_y\mu} - \frac{18\eta_x^4 nL^2\cC_3}{1-\lambda^2}\prt{1 + \frac{(4 + \eta_y\mu)\kappa^2}{\eta_y\mu}}}\mE\brk{\norm{\frac{1}{n}\sumn\nabla_x f_i(x_i^t, y_i^t)}^2}\\
			&\quad + \crk{\frac{4\kappa^2\cC_2}{n\eta_y\mu} + \cC_3\brk{2 + \frac{9\eta_y^2 L^2}{1-\lambda^2} + \frac{9\eta_x^2 L^2}{n(1-\lambda^2)}\prt{1 + \frac{(4 + \eta_y\mu)\kappa^2}{\eta_y\mu}}}}\frac{n\eta_x^3\sigma^2}{b}
		\end{aligned}
	\end{equation}

	To obtain the recursion for $\mE\cL^t$, it is sufficient to choose proper parameters to have 
	\begin{subequations}
		\label{eq:dsgda_lya1_1}
		\begin{align}
			\frac{L^2}{n} + \frac{1 + \lambda^2}{2}\cC_1 + \frac{9 L^2\eta_y\cC_2}{\mu} + \prt{2 + \frac{9\eta_y L^2}{\mu}}\frac{9\eta_x^2\cC_3}{1-\lambda^2} &\leq \cC_1\label{eq:c1}\\
			\frac{L^2}{n} + \prt{1 - \frac{\eta_y\mu}{4}}\cC_2 + \frac{18\eta_x^2 L^2\cC_3}{1-\lambda^2} &\leq \cC_2\label{eq:c2}\\
			\frac{2\lambda^2\cC_1}{1-\lambda^2} + \frac{3+\lambda^2}{4}\cC_3&\leq \cC_3.\label{eq:c3}
		\end{align}
	\end{subequations}

	From \eqref{eq:c3}, it suffices to choose $\cC_3 = 8\cC_1/(1-\lambda^2)^2$. Substituting such a $\cC_3$ into \eqref{eq:c2} and \eqref{eq:c1} leads to 
	\begin{subequations}
		\label{eq:dsgda_lya1_2}
		\begin{align}
			\frac{L^2}{n}  + \frac{9 L^2\eta_y\cC_2}{\mu} &\leq \brk{\frac{1-\lambda^2}{2} - \prt{2 + \frac{9\eta_y L^2}{\mu}}\frac{72\eta_x^2 L^2}{(1-\lambda^2)^3} }\cC_1\label{eq:c13}\\
			\frac{L^2}{n} + \frac{144\eta_x^2 L^2\cC_1}{(1-\lambda^2)^3} &\leq \frac{\eta_y\mu\cC_2}{4}.\label{eq:c23}
		\end{align}
	\end{subequations}

	Letting
	\begin{align*}
		\eta_y\leq \frac{1}{9\mu},\text{ and } \eta_x\leq\min\crk{\frac{(1-\lambda^2)^2}{12\sqrt{6}\kappa L}, \frac{(1-\lambda^2)^2}{120\sqrt{3}\kappa L}},
	\end{align*}
	it is sufficient to choose
	\begin{align}
		\label{eq:c1c2}
		\cC_1= \frac{300\kappa^2 L^2}{n(1-\lambda^2)},\; \cC_2 = \frac{8L^2}{n\eta_y\mu},\; \cC_3 = \frac{2400 \kappa^2 L^2}{n(1-\lambda^2)^3},
	\end{align}
	for \eqref{eq:dsgda_lya1_2} to hold.

	Letting $\eta_y\leq \sqrt{1-\lambda^2}/(3L)$ and $\eta_x\leq \min\crki{1/(12\kappa L), \eta_y/(8\sqrt{6}\kappa^2), (\eta_y\mu)^{1/4}(1-\lambda^2)/(50\kappa L)}$, the inequality \eqref{eq:dsgda_lya1_0} becomes 
	\begin{align*}
		\mE\cL^{t + 1}
		&\leq \mE\cL^t - \frac{\eta_x}{2}\mE\brk{\norm{\nabla \Phi(\bar{x}^2)}^2} + \frac{\eta_x^2 L\kappa \sigma^2}{nb} + \frac{24\eta_x\eta_y\kappa L\sigma^2}{b} \\
		&\quad + \brk{\frac{1}{n\eta_y^2\mu^2} + \frac{225}{(1-\lambda^2)^3} + \frac{54\eta_x^2\kappa^2 L^2}{n(1-\lambda^2)\eta_y\mu} }\frac{32\eta_x^3\kappa^2 L^2\sigma^2}{b},
	\end{align*}
    which completes the proof.

\end{document}